\documentclass[a4paper,12pt]{article}
\usepackage{amsmath, amssymb, latexsym, amscd, amsthm,amsfonts,amstext}
\usepackage[mathscr]{eucal}
\usepackage{graphicx}
\usepackage{subfig}
\usepackage{fancyhdr}
\usepackage{pgfplotstable}
\usepackage{tikz}
\usepackage[skins]{tcolorbox}

\numberwithin{equation}{section}
\input{tcilatex}
\usetikzlibrary{decorations.markings}
\usetikzlibrary{patterns}
\tikzset{onearrow/.style={decoration={markings, mark= at position 1.0 with {\arrow{#1}} , },postaction={decorate} }}

\begin{document}

\title{ Convexification of Restricted Dirichlet-to-Neumann Map }
\author{Michael V. Klibanov \and Department of Mathematics and Statistics,
University of North \and Carolina at Charlotte, Charlotte, NC 28223, USA
\and \texttt{mklibanv@uncc.edu} }
\date{}
\maketitle

\begin{abstract}
By our definition, \textquotedblleft restricted Dirichlet-to-Neumann map"
(DN) means that the Dirichlet and Neumann boundary data for a Coefficient
Inverse Problem (CIP) are generated by a point source running along an
interval of a straight line. On the other hand, the conventional DN data can
be generated, at least sometimes, by a point source running along a
hypersurface. CIPs with the restricted DN data are non-overdetermined in the 
$n-$D case with $n\geq 2.$ We develop, in a unified way, a \emph{general and
a radically new} numerical concept for CIPs with restricted DN data for a
broad class of PDEs of the second order, such as, e.g. elliptic, parabolic
and hyperbolic ones.\ Namely, using Carleman Weight Functions, we construct
globally convergent numerical methods. H\"{o}lder stability and uniqueness
are also proved. The price we pay for these features is a well acceptable
one in the Numerical Analysis: we truncate a certain Fourier-like series
with respect to some functions depending only on the position of that point
source. At least three applications are: imaging of land mines, crosswell
imaging and electrical impedance tomography.
\end{abstract}

\graphicspath{
{FIGURES/}
 {pics/}}



\bigskip

As of August 6, 2017, the corresponding paper is available online of Journal of Inverse and Ill-Posed Problems, DOI: 10.1515/jiip-2017-0067, \\
https://www.degruyter.com/printahead/j/jiip

\section{Introduction}

\label{sec:1}

The conventional Dirichlet-to-Neumann map (DN) data for a Coefficient
Inverse Problem (CIP) can be generated, at least sometimes, by the point
source running along a hypersurface, see pages 10-14 in \cite{KT} for DN and 
\cite{Hyv} for the Neumann-to-Dirichlet map data. We define
\textquotedblleft restricted DN data" \ for a CIP as the ones in which
Dirichlet and Neumann boundary data are generated by a point source running
along an interval of a straight line. These data are non-overdetermined in
the $n-$D case with $n\geq 2.$

We present in this paper a \emph{general and a radically new }concept of
constructing globally convergent numerical methods for CIPs with restricted
DN data. This concept also covers both H\"{o}lder stability and uniqueness
results for the CIPs we consider. Our construction is independent on a
specific PDE operator: it is the same for those PDEs of the second order,
which admit Carleman estimates. In particular, it works for three main types
of PDEs of the second order: elliptic, parabolic and hyperbolic ones. The
Dirichlet and Neumann data in elliptic and parabolic cases can be given on a
part of the boundary.

The price we pay for our concept is a well acceptable one in the Numerical
Analysis: we truncate a Fourier-like series with respect to a certain
orthonormal basis in the $L_{2}$ space of functions depending only on the
position of that point source. Next, to find spatially dependent
coefficients of that truncated series, we construct a weighted globally
strictly convex Tikhonov-like functional with the Carleman Weight Function
(CWF) in it. This is the function, which is involved in the Carleman
estimate for the corresponding PDE operator. Also, we establish the global
convergence of the gradient projection method to the exact solution under
the natural condition that the noise in the data tends to zero. As some
applications, we mention detection and identification of land mines,
crosswell imaging and electrical impedance tomography.

The construction of weighted strictly convex Tikhonov-like functionals with
CWFs in them was started by the author in 1997 \cite{Klib97} with the
recently renewed interest in \cite{BKconv,KNT,KK}.\ However, all these works
consider only CIPs with a single measurement data, as opposed to many
measurements of the current paper. In \cite{Kl} this technique was applied,
for the first time, to ill-posed Cauchy problems for a class of quasilinear
PDEs of the second order. The idea of \cite{Kl} was further explored in \cite%
{KB}. Numerical results can be found in \cite{KB,KNT}.

As to the DN data, a very substantial number of works have been published.
Since this paper is not a survey of DN, we refer to only a very few of them,
for brevity, and the reader can find other references in these publications.
Global uniqueness theorems for the elliptic case, i.e. for the Calderon
problem, were obtained in \cite{Nach,Nov1,SU}. Some reconstruction
procedures can be found in \cite{Lakh,Nach,Nov1,Nov3}. In the reconstruction
procedure of \cite{Nov3}, a certain infinite matrix is truncated, which is
philosophically close to our truncation of that Fourier-like series. We
refer to \cite{Alberti,Harrah,Jin} for numerical studies of DN. In \cite%
{Belishev} and \cite{Kabanikhin} reconstruction procedures for DN for
hyperbolic PDEs were developed, and they were computationally tested in \cite%
{Belishev1} and \cite{Kabanikhin,Kabanikhin1} respectively.

We point out that since our goal here is to present a new numerical concept,
for brevity, we are not concerned below with some issues related to
solutions of forward problems, since they can be discussed in later
publications. These issues are: the minimal smoothness assumptions,
existence and uniqueness of the solutions of the forward problems under
considerations, the positivity of those solutions and also the continuous
differentiability of those solutions with respect to the position of the
point source, see Conditions 1-3 in section 2.2. We just assume below that
these properties hold.

In sections 2-4 we present our concept for the case of a general PDE of the
second order, for which a Carleman estimate is valid. Next, we specify this
concept for elliptic, parabolic and hyperbolic PDEs in sections 5, 6 and 7
respectively. In particular, we outline in section 5 applications to
detection and identification of land mines, crosswell imaging and electrical
impedance tomography. Finally, we present in section 8 some thoughts about
numerical studies.

\section{A CIP With the Restricted DN Data}

\label{sec:2}

\subsection{The Carleman estimate}

\label{sec:2.1}

Below all functions are real valued, unless stated otherwise. The material
of section 2.1 is a somewhat modified material of section 2.1.2 of \cite{KT}%
. Below $x=\left( x_{1},...,x_{n}\right) \in \mathbb{R}^{n}$. Also, below $%
\alpha =\left( \alpha _{1},...,\alpha _{n}\right) $ is the multiindex with
integer coordinates $\alpha _{i}\geq 0$ and with $\left\vert \alpha
\right\vert =\alpha _{1}+...+\alpha _{n}.$ Consider a general Partial
Differential Operator of the second order%
\begin{eqnarray}
A\left( x,u\right) &=&\sum\limits_{\left\vert \alpha \right\vert \leq
2}a_{\alpha }\left( x\right) D_{x}^{\alpha }u=A_{0}\left( x,u\right)
+A_{1}\left( x,u\right) ,x\in \mathbb{R}^{n},  \label{2.1} \\
A_{0}\left( x,u\right) &=&\sum\limits_{\left\vert \alpha \right\vert
=2}a_{\alpha }\left( x\right) D_{x}^{\alpha }u,\text{ }A_{1}\left(
x,u\right) =\sum\limits_{\left\vert \alpha \right\vert =1}a_{\alpha }\left(
x\right) D_{x}^{\alpha }u+a_{0}\left( x\right) u.  \label{2.2}
\end{eqnarray}%
Thus, $A_{0}\left( x,u\right) $ is the principal part of the operator $%
A\left( x,u\right) $ and the operator $A_{1}\left( x,u\right) $ contains
lower order terms. Let $\Omega \subset \mathbb{R}^{n}$ be a bounded domain
with a piecewise smooth boundary. Let $Z>0$ be a given number. We assume
that coefficients%
\begin{eqnarray}
\text{ }a_{\alpha }\left( x\right) &=&\widehat{a}_{\alpha }=const.\text{ for 
}x\notin \Omega \text{ and for all }\alpha \text{ with }\left\vert \alpha
\right\vert \leq 2,  \label{2.30} \\
a_{\alpha } &\in &C^{1}\left( \mathbb{R}^{n}\right) \text{ for }\left\vert
\alpha \right\vert \leq 2,  \label{2.31} \\
\left\Vert a_{\alpha }\right\Vert _{C^{1}\left( \overline{\Omega }\right) }
&\leq &Z\text{ for }\left\vert \alpha \right\vert \leq 1.  \label{2.4}
\end{eqnarray}

Let $\Gamma \in C^{2},\Gamma \subseteq \partial \Omega $ be a part of the
boundary of the domain $\Omega .$ We assume that any part of $\Gamma $ is
not a characteristic surface of the operator $A_{0}\left( x,u\right) .$ Let
the function $\xi \in C^{\infty }\left( \overline{\Omega }\right) $ and $%
\left\vert \nabla \xi \right\vert \neq 0$ in $\overline{\Omega }.$ For a
number $d>0$ denote 
\begin{equation}
\xi _{d}=\left\{ x\in \Omega :\xi \left( x\right) =d\right\} ,\Omega
_{d}=\left\{ x\in \Omega :\xi \left( x\right) >d\right\} .  \label{2.6}
\end{equation}%
We assume below that $\Omega _{d}\neq \varnothing $ and that $\left( 
\overline{\Omega }_{d}\cap \partial \Omega \right) =\Gamma _{d}\subseteq
\Gamma .$ Hence, 
\begin{equation}
\Gamma _{d}=\left\{ x\in \Gamma :\xi \left( x\right) >d\right\} .
\label{2.7}
\end{equation}%
Hence, the boundary of the domain $\Omega _{d}$ consists of two parts,%
\begin{equation}
\partial \Omega _{d}=\partial _{1}\Omega _{d}\cup \partial _{2}\Omega _{d},%
\text{ }\partial _{1}\Omega _{d}=\xi _{d},\partial _{2}\Omega _{d}=\Gamma
_{d}.  \label{2.8}
\end{equation}

We assume below that $\partial \Omega _{d}$ is piecewise smooth. Below $%
C_{1}=C_{1}\left( A_{0},\Omega _{d}\right) >0$ denotes different constants
depending only on the operator $A_{0}$ and the domain $\Omega .$ Let $%
\lambda >1$ be a large parameter. Consider the function $\varphi _{\lambda
}\left( x\right) ,$%
\begin{equation}
\varphi _{\lambda }\left( x\right) =\exp \left( \lambda \xi \left( x\right)
\right) .  \label{2}
\end{equation}%
It follows from (\ref{2.6})-(\ref{2.8}) that 
\begin{equation}
\min_{\overline{\Omega _{d}}}\varphi _{\lambda }\left( x\right) =\varphi
_{\lambda }\left( x\right) \mid _{\xi _{d}}=e^{\lambda d},  \label{2.80}
\end{equation}%
\begin{equation}
m=\max_{\overline{\Omega _{d}}}\xi \left( x\right) \Rightarrow \max_{%
\overline{\Omega _{d}}}\varphi _{\lambda }\left( x\right) =e^{\lambda m}.
\label{2.81}
\end{equation}

\textbf{Definition 2.1}. \emph{We say that the operator }$A_{0}$\emph{\ with
its coefficients }$a_{\alpha }\left( x\right) $ \emph{satisfying conditions} 
\emph{(\ref{2.2}), (\ref{2.31})} \emph{admits the pointwise\ Carleman
estimate in the domain }$\Omega _{d}$\emph{\ with the CWF }$\varphi
_{\lambda }\left( x\right) $\emph{\ if there exist constants }$\lambda
_{0}=\lambda _{0}\left( A_{0},\Omega _{d}\right) >1,C_{1}=C_{1}\left(
A_{0},\Omega _{d}\right) >0,$\emph{\ depending only on listed parameters,
such that the following estimates hold}%
\begin{eqnarray}
\left( A_{0}u\right) ^{2}\varphi _{\lambda }^{2}\left( x\right) &\geq
&C_{1}\lambda \left( \nabla u\right) ^{2}\varphi _{\lambda }^{2}\left(
x\right) +C_{1}\lambda ^{3}u^{2}\varphi _{\lambda }^{2}\left( x\right) +%
\func{div}U,  \label{2.90} \\
\left\vert U\left( x\right) \right\vert &\leq &C_{1}\lambda ^{3}\left[
\left( \nabla u\right) ^{2}+u^{2}\right] \varphi _{\lambda }^{2}\left(
x\right) ,  \label{2.91} \\
\forall \lambda &\geq &\lambda _{0},\forall x\in \overline{\Omega _{d}}%
,\forall u\in C^{2}\left( \overline{\Omega _{d}}\right) .  \label{2.92}
\end{eqnarray}

\subsection{Statement of the problem}

\label{sec:2.2}

Denote $\overline{x}=\left( x_{2},...,x_{n}\right) \in \mathbb{R}^{n-1}.$
Below $\overline{x}^{0}\in \mathbb{R}^{n-1}$ is a fixed point of $\mathbb{R}%
^{n-1}$ and $x_{0}\in \left[ 0,1\right] $ is a varying parameter. Consider
an interval $I$ of a straight line such that%
\begin{eqnarray}
I &=&\left\{ x=\left( x_{0},\overline{x}^{0}\right) :x_{0}\in \left(
0,1\right) \right\} ,  \label{2.10} \\
I\cap \overline{\Omega } &=&\varnothing .  \label{2.11}
\end{eqnarray}%
Consider the following equation%
\begin{equation}
A\left( u\right) =-\delta \left( x_{1}-x_{0},\overline{x}-\overline{x}%
^{0}\right) ,x\in \mathbb{R}^{n},\forall x_{0}\in \left[ 0,1\right] ,
\label{2.12}
\end{equation}%
where $u=u\left( x,x_{0}\right) $ is a distribution with respect to $x.$
Since we do not impose any condition at the infinity on the distribution $u$%
, equation (\ref{2.12}) might have many solutions or even none. Suppose that
it has a solution, which we still denote as $u\left( x,x_{0}\right) .$ We
assume that the following conditions are valid for this solution:

\textbf{Condition 1}. For each $x_{0}\in \left[ 0,1\right] $ the function $%
u\left( x,x_{0}\right) \in C^{2}\left( \overline{\Omega }\right) .$

\textbf{Condition 2. }For each $x\in \overline{\Omega },$ the functions $%
D_{x}^{\alpha }u\left( x,x_{0}\right) $, are differentiable with respect to $%
x_{0}\in \left( 0,1\right) $ and functions $\partial
_{x_{0}}^{k}D_{x}^{\alpha }u\left( x,x_{0}\right) \in C\left( \overline{%
\Omega }\times \left[ 0,1\right] \right) $ for $k=0,1;$ $\left\vert \alpha
\right\vert \leq 2.$

\textbf{Condition 3. }$u\left( x,x_{0}\right) \geq \beta =const.>0,\forall
\left( x,x_{0}\right) \in \overline{\Omega }\times \left[ 0,1\right] ,$ see
Remark 2.1.

\textbf{Condition 4. }The following Dirichlet and Neumann boundary
conditions hold for the function $u\left( x,x_{0}\right) :$%
\begin{equation}
u\left( x,x_{0}\right) \mid _{x\in \Gamma ,x_{0}\in \left[ 0,1\right]
}=g_{0}\left( x,x_{0}\right) ,\partial _{n}u\left( x,x_{0}\right) \mid
_{x\in \Gamma ,x_{0}\in \left[ 0,1\right] }=g_{1}\left( x,x_{0}\right) ,
\label{2.130}
\end{equation}%
where $g_{0}\left( x,x_{0}\right) $ and $g_{1}\left( x,x_{0}\right) $ are
two given functions of $\left( x,x_{0}\right) \in \Gamma \times \left[ 0,1%
\right] .$

We call the Dirichlet and Neumann boundary data (\ref{2.130})\emph{\ }%
\textquotedblleft restricted DN data".

\textbf{Coefficient Inverse Problem 1 (CIP1). }\emph{Suppose that for each
value }$x_{0}\in \left[ 0,1\right] $ \emph{of the} \emph{parameter }$x_{0}$ 
\emph{there exists a distribution }$u\left( x,x_{0}\right) $\emph{\
satisfying equation (\ref{2.12}) and Conditions 1-4. Determine the
coefficient }$a_{0}\left( x\right) $\emph{\ in (\ref{2.2}) from functions }$%
g_{0}\left( x,x_{0}\right) $\emph{\ and }$g_{1}\left( x,x_{0}\right) $\emph{%
\ in (\ref{2.130}).}

\textbf{Remark 2.1.} \emph{Thus, (\ref{2.12}) and (\ref{2.130}) mean that
the source }$\left( x_{0},\overline{x}^{0}\right) $ \emph{runs along the
interval }$I$. \emph{In the cases of elliptic and parabolic PDEs Condition 3
can often be established via the maximum principle \cite{Fr,G}. }

Sometimes it is hard to prove the validity of Conditions 1-3 in the case
when the fundamental solution (\ref{2.12}) of the operator $A$ is
considered. Hence, we formulate now the second CIP with restricted DN data.
Let $\varepsilon >0$ be a sufficiently small number. Let the functions $f\in
C^{\infty }\left( \mathbb{R}\right) $ and $\chi \left( \overline{x}\right)
\in C^{\infty }\left( \mathbb{R}^{n-1}\right) $ be such that $f\left(
0\right) \chi \left( 0\right) \neq 0$ and also $f\left( z\right) =0$ for $%
\left\vert z\right\vert >\varepsilon $ as well as $\chi \left( y\right) =0$
for $y\in \left\{ \left\vert y\right\vert >\varepsilon \right\} .$ Let $%
I_{\varepsilon }=\left\{ x\in \mathbb{R}^{n}:dist\left( x,I\right)
<\varepsilon \right\} ,$ where $dist\left( x,I\right) $ is the Hausdorff
distance between the point $x$ and the interval $I$. Let $G\subset \mathbb{R}%
^{n}$ be a bounded domain with its boundary $\partial G\in C^{1}$ and such
that $\Omega \subset G,\partial \Omega \cap \partial G=\varnothing .$ We
assume that $I_{\varepsilon }\subset \left( G\diagdown \overline{\Omega }%
\right) .$

We now replace (\ref{2.12}) with%
\begin{equation}
A\left( \widetilde{u}\right) =f\left( x_{1}-x_{0}\right) \chi \left( 
\overline{x}-\overline{x}^{0}\right) ,\forall x_{0}\in \left[ 0,1\right] ,
\label{2.14}
\end{equation}%
\begin{equation}
\widetilde{u}\mid _{x\in \partial G}=0,\forall x_{0}\in \left[ 0,1\right] .
\label{2.15}
\end{equation}

\textbf{Coefficient Inverse Problem 2 (CIP2).} \emph{Assume that the
function }$\widetilde{u}\left( x,x_{0}\right) $\emph{\ satisfies Conditions
1-4, (\ref{2.14}) and (\ref{2.15}). Determine the coefficient }$a_{0}\left(
x\right) $\emph{\ in (\ref{2.2}) from functions }$g_{0}\left( x,x_{0}\right) 
$\emph{\ and }$g_{1}\left( x,x_{0}\right) $\emph{\ in (\ref{2.130}). }

Both CIP1 and CIP2 are non overdetermined.\ Indeed the number $n$ of free
variables in the data (\ref{2.130}) coincides with the number of free
variables in the unknown coefficient. Since our method of the numerical
solution of CIP2 is exactly the same as the one of CIP1, we consider below
CIP1 in most cases.

\subsection{A special orthonormal basis in $L_{2}\left( 0,1\right) $}

\label{sec:2.3}

We need to construct such an orthonormal basis in the space $L_{2}\left(
0,1\right) $ of functions depending on $x_{0}$ that the first derivative
with respect to $x_{0}$ of any element of this basis is not identically
zero. In addition, this derivative should be a linear combination of a
finite number of elements of this basis. Neither the basis of trigonometric
functions nor the basis of standard orthonormal polynomials are not suitable
for this goal. Therefore, we construct a new basis. Our basis is similar
with Laguerre functions, which, however, form an orthonormal basis in $%
L_{2}\left( 0,\infty \right) $ rather than in $L_{2}\left( 0,1\right) .$

For $x_{0}\in \left( 0,1\right) ,$ consider the set of functions $\left\{
x_{0}^{k}e^{x_{0}}\right\} _{k=0}^{\infty }.$ Clearly these functions are
linearly independent and form a complete set in $L_{2}\left( 0,1\right) .$
We apply the classical Gram-Schmidt Orthonormalization procedure to this
set. We start from $e^{x_{0}}.$ Then we take $x_{0}e^{x_{0}},$ then $%
x_{0}^{2}e^{x_{0}},$ etc. As a result, we obtain an orthonormal basis in $%
L_{2}\left( 0,1\right) ,$ which consists of functions $\left\{ P_{m}\left(
x_{0}\right) e^{x_{0}}\right\} _{m=0}^{\infty }=\left\{ \psi _{m}\left(
x_{0}\right) \right\} _{m=0}^{\infty },$ where $P_{m}\left( x_{0}\right) $
is a polynomial of the degree $m$. Denote $\left[ ,\right] $ the scalar
product in $L_{2}\left( 0,1\right) .$ Let $Q_{s}\left( x_{0}\right) $ be an
arbitrary polynomial of the degree $s\geq 0$. By the construction of
functions $\psi _{m}\left( x_{0}\right) ,$ there exists numbers $%
b_{j}=b_{j}\left( Q_{s}\right) $ such that 
\begin{equation}
Q_{s}\left( x_{0}\right) =\displaystyle\sum \limits_{j=0}^{s}b_{j}\left(
Q_{s}\right) P_{j}\left( x_{0}\right) .  \label{2.16}
\end{equation}

\textbf{Remark 2.2.} \emph{In the computational practice, one can use a
symbolic computations software, such as, e.g. Mathematica or Maple to figure
out a reasonable number of functions }$\psi _{m}\left( x_{0}\right) .$

\textbf{Theorem 2.1.} \emph{We have} 
\begin{equation}
a_{mk}=\left[ \psi _{k}^{\prime },\psi _{m}\right] =\left\{ 
\begin{array}{c}
1\text{ if }k=m, \\ 
0\text{ if }k<m.%
\end{array}%
\right.  \label{2.17}
\end{equation}%
\emph{Let }$N>1$\emph{\ be an integer. Consider the }$N\times N$\emph{\
matrix }$M_{N}=\left( a_{mk}\right) _{\left( k,m\right) =\left( 0,0\right)
}^{\left( N-1,N-1\right) }.$\emph{\ Then (\ref{2.17}) implies that }$\det
\left( M_{N}\right) =1,$ \emph{which means that there exists the inverse
matrix }$M_{N}^{-1}.$

\textbf{Proof}. We have $\psi _{k}^{\prime }\left( x_{0}\right) =P_{k}\left(
x_{0}\right) e^{x_{0}}+P_{k}^{\prime }\left( x_{0}\right) e^{x_{0}}=\psi
_{k}\left( x_{0}\right) +P_{k}^{\prime }\left( x_{0}\right) e^{x_{0}}.$
Since the degree of the polynomial $P_{k}^{\prime }\left( x_{0}\right) $ is
less than $k$, then (\ref{2.16}) implies that the function $P_{k}^{\prime
}\left( x_{0}\right) e^{x_{0}}$ is a linear combination of functions $\psi
_{j}\left( x_{0}\right) $ with $j\leq k-1.$ Hence, 
\begin{equation}
\psi _{k}^{\prime }\left( x_{0}\right) =\psi _{k}\left( x_{0}\right) +%
\displaystyle\sum \limits_{j=0}^{k-1}b_{jk}\psi _{j}\left( x_{0}\right) .
\label{2.18}
\end{equation}%
First, let $m=k.$ Since $\left[ \psi _{j},\psi _{m}\right] =0$ for $j<m,$
then (\ref{2.18}) implies that $\left[ \psi _{m}^{\prime }\left(
x_{0}\right) ,\psi _{m}\left( x_{0}\right) \right] =1.$ Consider now the
case $m>k.$ Then we obtain similarly from (\ref{2.18}) that

$\left[ \psi _{k}^{\prime }\left( x_{0}\right) ,\psi _{m}\left( x_{0}\right) %
\right] =0.$ Thus, (\ref{2.17}) is established. $\square $

\section{An Ill-Posed Problem for a Coupled System of Quasilinear PDEs}

\label{sec:3}

If we say below that a certain vector function belongs to a functional
space, then this means that each component of this function belongs to this
space. The norm of that vector function in that space is defined as the
square root of the sum of squares of norms of its components.

It follows from Condition 3 of section 2.2 that we can consider the function 
$v\left( x,x_{0}\right) =\ln u\left( x,x_{0}\right) $ for $x\in \overline{%
\Omega }.$ Substituting $u=e^{v}$ in (\ref{2.12}) for $x\in \Omega $ and
using (\ref{2.1})-(\ref{2.4}), (\ref{2.10}), (\ref{2.11}) and (\ref{2.130}),
we obtain%
\begin{equation}
A_{0}\left( x,v\right) +F_{1}\left( x,\nabla v\right) =-a_{0}\left( x\right)
,x\in \Omega ,x_{0}\in \left[ 0,1\right] ,  \label{3.1}
\end{equation}%
\begin{equation}
v\left( x,x_{0}\right) \mid _{x\in \Gamma ,x_{0}\in \left[ 0,1\right] }=%
\widetilde{g}_{0}\left( x,x_{0}\right) ,\partial _{n}v\left( x,x_{0}\right)
\mid _{x\in \Gamma ,x_{0}\in \left[ 0,1\right] }=\widetilde{g}_{1}\left(
x,x_{0}\right) ,  \label{3.2}
\end{equation}%
\begin{equation*}
\widetilde{g}_{0}\left( x,x_{0}\right) =\ln g_{0}\left( x,x_{0}\right) ,%
\widetilde{g}_{1}\left( x,x_{0}\right) =g_{1}\left( x,x_{0}\right)
/g_{0}\left( x,x_{0}\right) ,
\end{equation*}%
where the function $F_{1}\in C^{1}\left( \mathbb{R}^{2n}\right) ,$ and it is
quadratic with respect to derivatives $\partial _{x_{k}}v.$ Denote $%
v_{x_{0}}\left( x,x_{0}\right) =\partial _{x_{0}}v\left( x,x_{0}\right) .$
Differentiate both sides of (\ref{3.1}) with respect to $x_{0}.$ Since $%
\partial _{x_{0}}\left( a_{0}\left( x\right) \right) \equiv 0,$ then, using (%
\ref{3.2}), we obtain%
\begin{equation}
A_{0}\left( x,v_{x_{0}}\right) +F_{2}\left( x,\nabla v,\nabla
v_{x_{0}}\right) =0,x\in \Omega ,x_{0}\in \left[ 0,1\right] ,  \label{3.3}
\end{equation}%
\begin{equation}
v_{x_{0}}\left( x,x_{0}\right) \mid _{x\in \Gamma ,x_{0}\in \left[ 0,1\right]
}=\partial _{x_{0}}\widetilde{g}_{0}\left( x,x_{0}\right) ,\partial
_{n}v_{x_{0}}\left( x,x_{0}\right) \mid _{x\in \Gamma ,x_{0}\in \left[ 0,1%
\right] }=\partial _{x_{0}}\widetilde{g}_{1}\left( x,x_{0}\right) ,
\label{3.4}
\end{equation}%
where the function $F_{2}\in C^{1}\left( \mathbb{R}^{3n}\right) $ is
quadratic with respect to derivatives $\partial _{x_{k}}v,$ $\partial
_{x_{k}}v_{x_{0}}.$

It follows from Conditions 1-3 of section 2.2 that, for each $x\in \overline{%
\Omega },$ the function $v\left( x,x_{0}\right) $ can be represented as a
Fourier-like series with respect to the orthonormal basis $\left\{ \psi
_{m}\left( x_{0}\right) \right\} _{m=0}^{\infty }$. Coefficients of this
series depend on $x$. We, however, assume that the function $v\left(
x,x_{0}\right) $ can be represented as a truncated series, 
\begin{equation}
v\left( x,x_{0}\right) =\displaystyle\sum\limits_{k=0}^{N-1}v_{k}\left(
x\right) \psi _{k}\left( x_{0}\right) ,\forall x\in \overline{\Omega }%
,\forall x_{0}\in \left[ 0,1\right] ,  \label{3.5}
\end{equation}%
where coefficients $v_{k}\left( x\right) \in C^{2}\left( \overline{\Omega }%
\right) $ and $N\geq 2$ is an integer of ones choice. Substituting (\ref{3.5}%
) in (\ref{3.3}), we obtain 
\begin{equation}
\displaystyle\sum\limits_{k=0}^{N-1}A_{0}\left( x,v_{k}\right) \psi
_{k}^{\prime }\left( x_{0}\right) +F_{2}\left( x,\displaystyle%
\sum\limits_{m=0}^{N-1}\nabla v_{m}\left( x\right) \psi _{m}\left(
x_{0}\right) ,\displaystyle\sum\limits_{k=0}^{N-1}\nabla v_{k}\left(
x\right) \psi _{k}^{\prime }\left( x_{0}\right) \right) =0,  \label{3.6}
\end{equation}%
where $x\in \Omega ,x_{0}\in \left[ 0,1\right] .$ Let the integer $m\in %
\left[ 0,N-1\right] .$ Multiply both sides of (\ref{3.6}) by the function $%
\psi _{m}\left( x_{0}\right) $ and integrate with respect to $x_{0}\in
\left( 0,1\right) .$ We obtain%
\begin{equation*}
\displaystyle\sum\limits_{k=0}^{N-1}A_{0}\left( x,v_{k}\right) \left[ \psi
_{k}^{\prime }\left( x_{0}\right) ,\psi _{m}\left( x_{0}\right) \right]
\end{equation*}%
\begin{equation}
=-\left[ F_{2}\left( x,\displaystyle\sum\limits_{k=0}^{N-1}\nabla
v_{k}\left( x\right) \psi _{k}\left( x_{0}\right) ,\displaystyle%
\sum\limits_{k=0}^{N-1}\nabla v_{k}\left( x\right) \psi _{k}^{\prime }\left(
x_{0}\right) \right) ,\psi _{m}\left( x_{0}\right) \right] ,  \label{3.7}
\end{equation}%
where $x\in \Omega ,m\in \left[ 0,N-1\right] .$ Denote 
\begin{equation}
V\left( x\right) =\left( v_{0}\left( x\right) ,...,v_{N-1}\left( x\right)
\right) ^{T},A_{0}\left( x,V\right) =\left( A_{0}\left( x,v_{0}\right)
,...,A_{0}\left( x,v_{N-1}\right) \right) ^{T}.  \label{3.8}
\end{equation}%
Also, let $F\left( x,\nabla V\right) =\left( F_{2,0}\left( x,\nabla V\right)
,...,F_{2,N-1}\left( x,\nabla V\right) \right) ^{T}$ be the vector of right
hand sides of equations (\ref{3.7}). Then (\ref{3.7}) can be rewritten as%
\begin{equation}
M_{N}A_{0}\left( x,V\right) =F\left( x,\nabla V\right) ,  \label{3.9}
\end{equation}%
where $M_{N}$ is the matrix of Theorem 2.1. Applying Theorem 2.1 to (\ref%
{3.9}) and denoting $P\left( x,\nabla V\right) =M_{N}^{-1}F\left( x,\nabla
V\right) $, we obtain%
\begin{equation}
A_{0}\left( x,V\right) -P\left( x,\nabla V\right) =0,  \label{3.10}
\end{equation}%
\begin{equation}
V\mid _{\Gamma }=p_{0}\left( x\right) ,\partial _{n}V\mid _{\Gamma
}=p_{1}\left( x\right) ,  \label{3.11}
\end{equation}%
where vector functions $p_{0}\left( x\right) $ and $p_{1}\left( x\right) $
are obtained from functions $\partial _{x_{0}}\widetilde{g}_{0}\left(
x,x_{0}\right) $ and $\partial _{x_{0}}\widetilde{g}_{1}\left(
x,x_{0}\right) $ of (\ref{3.4}) in an obviously similar manner, the $N-$%
Dimensional vector function $P\in C^{1}\left( \mathbb{R}^{s_{1}}\right)
,s_{1}=n\left( N+1\right) ,$ and each component of $P$ is a quadratic
function of the first derivatives $\partial _{x_{k}}v_{i}\left( x\right) ,$
where $k=1,...,n$ and $i=0,...,N-1.$

Equalities (\ref{3.10}), (\ref{3.11}) form an ill-posed problem for the
coupled system of quasilinear equations. A similar problem was considered in 
\cite{KB,Kl} for the case of a single quasilinear PDE. Thus, we proceed
below similarly with \cite{KB,Kl}. It follows from (\ref{3.1}), (\ref{3.5})
and (\ref{3.8}) that, given the vector function $V\left( x\right) $, we can
find the unknown coefficient $a_{0}\left( x\right) .$ However, since only $u$
and $\nabla u$ are involved in the Carleman estimate (\ref{2.90}), while the
second derivatives $u_{x_{i}x_{j}}$ are not involved, we formulate all
theorems below in terms of the vector functions $V,\nabla V$ rather than in
terms of the unknown coefficient $a_{0}\left( x\right) .$ At the same time,
it is well known that in the case of parabolic and elliptic operators
(unlike hyperbolic ones) derivatives involved in their principal parts can
be incorporated in corresponding Carleman estimates, see, e.g. theorem 2.5
in \cite{Ksurvey}. Hence, the H\"{o}lder stability result of Theorem 3.1 as
well as the global convergence result (Theorem 4.4 below) can be
reformulated in terms of $a_{0}\left( x\right) $ in these cases. We are not
doing this here for brevity.

\textbf{Theorem 3.1} (H\"{o}lder stability and uniqueness). \emph{Suppose
that there exist two vector functions }$V^{\left( 1\right) },V^{\left(
2\right) }\in C^{2}\left( \overline{\Omega }\right) $\emph{\ satisfying
equation (\ref{3.10}) and with two pairs of boundary conditions (\ref{3.11}%
), }$V^{\left( i\right) }\mid _{\Gamma }=p_{0}^{\left( i\right) }\left(
x\right) $\emph{\ and }$\partial _{n}V^{\left( i\right) }\mid _{\Gamma
}=p_{1}^{\left( i\right) }\left( x\right) ,i=1,2.$\emph{\ Let }$K>0$\emph{\
be such a number that }$\left\Vert V^{\left( i\right) }\right\Vert
_{C^{1}\left( \overline{\Omega }\right) }\leq K$\emph{. Let }$Z>0$\emph{\ be
the number defined in (\ref{2.4}). Let }$\sigma \in \left( 0,1\right) $\emph{%
\ be the level of the error in the data (\ref{3.11}), i.e.}%
\begin{equation}
\left\Vert p_{0}^{\left( 1\right) }-p_{0}^{\left( 2\right) }\right\Vert
_{H^{1}\left( \Gamma \right) }\leq \sigma ,\left\Vert p_{1}^{\left( 1\right)
}-p_{1}^{\left( 2\right) }\right\Vert _{L_{2}\left( \Gamma \right) }\leq
\sigma .  \label{3.12}
\end{equation}%
\emph{Choose a number }$c>0$\emph{\ such that }$\Omega _{d+c}\neq
\varnothing .$\emph{\ Then there exists a sufficiently small constant }$%
\sigma _{0}=\sigma _{0}\left( \Omega ,K,Z,\xi ,m,c\right) \in \left(
0,1\right) $\emph{\ and a constant }$C_{2}=C_{2}\left( \Omega ,K,Z,\xi
,m,c\right) >0,$\emph{\ both depending only on listed parameters, such that
for all }$\sigma \in \left( 0,\sigma _{0}\right) $\emph{\ the following H%
\"{o}lder stability estimate is valid} 
\begin{equation}
\left\Vert V^{\left( 1\right) }-V^{\left( 2\right) }\right\Vert
_{H^{1}\left( \Omega _{d+c}\right) }\leq C_{2}\left( 1+\left\Vert V^{\left(
1\right) }-V^{\left( 2\right) }\right\Vert _{H^{2}\left( \Omega \right)
}\right) \sigma ^{\rho },\rho =c/\left( m+c\right) .  \label{3.13}
\end{equation}%
\emph{In particular, if }$p_{0}^{\left( 1\right) }=p_{0}^{\left( 2\right) }$%
\emph{\ and }$p_{1}^{\left( 1\right) }=p_{1}^{\left( 2\right) },$\emph{\
i.e. if }$\sigma =0,$\emph{\ then }$V^{\left( 1\right) }\left( x\right)
=V^{\left( 2\right) }\left( x\right) $\emph{\ in }$\Omega _{d},$ \emph{which
means that uniqueness of the problem (\ref{3.10}), (\ref{3.11}) holds in the
domain }$\Omega _{d}.$

\textbf{Proof}. Uniqueness in the domain $\Omega _{d}$ follows from (\ref%
{3.13}) immediately.\emph{\ }In this proof, $C_{2}=C_{2}\left( \Omega
,K,Z,\xi ,m,c\right) >0$ denotes different positive constants$.$ Consider
the set of vector functions $Y=Y\left( K\right) =$ $\left\{ V\in C^{2}\left( 
\overline{\Omega }\right) :\left\Vert V\right\Vert _{C^{1}\left( \overline{%
\Omega }\right) }\leq K\right\} .$\ Denote $\widetilde{V}\left( x\right)
=V^{\left( 1\right) }\left( x\right) -V^{\left( 2\right) }\left( x\right) .$
Then $\widetilde{V}\left( x\right) =\left( \widetilde{v}_{0}\left( x\right)
,...,\widetilde{v}_{N-1}\left( x\right) \right) ^{T}.$ Since each component
of the vector function $P\left( x,\nabla V\right) $ is a quadratic function
with respect to the first derivatives $\partial _{x_{k}}v_{i}\left( x\right) 
$, then 
\begin{equation}
P\left( x,\nabla V^{\left( 1\right) }\right) -P\left( x,\nabla V^{\left(
2\right) }\right) =\widehat{P}\left( x,\nabla V^{\left( 1\right) },\nabla
V^{\left( 2\right) }\right) \nabla \widetilde{V}\left( x\right) ,
\label{3.14}
\end{equation}%
where the matrix $\widehat{P}\left( x,\nabla V^{\left( 1\right) },\nabla
V^{\left( 2\right) }\right) $ is such that 
\begin{equation}
\max_{V^{\left( 1\right) },V^{\left( 2\right) }\in Y}\left\Vert \widehat{P}%
\left( x,\nabla V^{\left( 1\right) },\nabla V^{\left( 2\right) }\right)
\right\Vert _{C\left( \overline{\Omega }\right) }\leq C_{2}.  \label{3.15}
\end{equation}

We obtain from (\ref{3.10}), (\ref{3.11}), (\ref{3.14}) and (\ref{3.15}) 
\begin{equation}
\left\vert A_{0}\left( x,\widetilde{v}_{i}\right) \right\vert \leq
C_{2}\left\vert \nabla \widetilde{V}\left( x\right) \right\vert ,\forall
x\in \Omega ,i=0,...,N-1,  \label{3.16}
\end{equation}%
\begin{equation}
\widetilde{v}_{i}\mid _{\Gamma }=\widetilde{p}_{0,i}\left( x\right)
,\partial _{n}\widetilde{v}_{i}\mid _{\Gamma }=\widetilde{p}_{1,i}\left(
x\right) ,i=0,...,N-1,  \label{3.17}
\end{equation}%
where $\widetilde{p}_{k,i}\left( x\right) =\left( p_{k,i}^{\left( 1\right)
}-p_{k,i}^{\left( 2\right) }\right) \left( x\right) ,k=0,1.$ Square both
sides of (\ref{3.16}), sum up with respect to $i=0,...,N-1,$ multiply by the
function $\varphi _{\lambda }^{2}\left( x\right) $ defined in (\ref{2}),
integrate over the domain $\Omega _{d}$, and then apply (\ref{2.8}), (\ref%
{2.90})-(\ref{2.92}) as well as the Gauss' formula. Also, use (\ref{2.80}), (%
\ref{2.81}) and (\ref{3.12}). Since $\left\Vert \widetilde{V}\right\Vert
_{L_{2}\left( \xi _{d}\right) },\left\Vert \nabla \widetilde{V}\right\Vert
_{L_{2}\left( \xi _{d}\right) }\leq C_{1}\left\Vert \widetilde{V}\right\Vert
_{H^{2}\left( \Omega \right) },$ then we obtain for $\lambda \geq \lambda
_{0}$%
\begin{equation}
C_{1}\lambda \displaystyle\int\limits_{\Omega _{d}}\left\vert \nabla 
\widetilde{V}\left( x\right) \right\vert ^{2}\varphi _{\lambda
}^{2}dx+C_{1}\lambda ^{3}\displaystyle\int\limits_{\Omega _{d}}\left\vert 
\widetilde{V}\left( x\right) \right\vert ^{2}\varphi _{\lambda }^{2}dx
\label{3.18}
\end{equation}%
\begin{equation*}
\leq C_{1}\lambda ^{3}e^{2\lambda m}\sigma ^{2}+C_{1}\lambda ^{3}e^{2\lambda
d}\left\Vert \widetilde{V}\right\Vert _{H^{2}\left( \Omega \right)
}^{2}+C_{2}\displaystyle\int\limits_{\Omega _{d}}\left\vert \nabla 
\widetilde{V}\left( x\right) \right\vert ^{2}\varphi _{\lambda }^{2}dx.
\end{equation*}%
Choose $\lambda _{1}=\lambda _{1}\left( \Omega ,K,Z,c\right) >\max \left(
\lambda _{0},1\right) $ so large that $C_{2}<C_{1}\lambda _{1}/2.$ Then we
obtain from (\ref{3.18}) for $\lambda \geq \lambda _{1}$%
\begin{equation}
\lambda \displaystyle\int\limits_{\Omega _{d}}\left\vert \nabla \widetilde{V}%
\left( x\right) \right\vert ^{2}\varphi _{\lambda }^{2}dx+\lambda ^{3}%
\displaystyle\int\limits_{\Omega _{d}}\left\vert \widetilde{V}\left(
x\right) \right\vert ^{2}\varphi _{\lambda }^{2}dx\leq C_{1}\lambda
^{3}e^{2\lambda m}\sigma ^{2}+C_{1}\lambda ^{3}e^{2\lambda d}\left\Vert 
\widetilde{V}\right\Vert _{H^{2}\left( \Omega \right) }^{2}.  \label{3.19}
\end{equation}%
Since $\Omega _{d+c}\subset \Omega _{d}$, $\Omega _{d+c}\neq \varnothing $
and also since 
\begin{equation}
\varphi _{\lambda }^{2}\left( x\right) >e^{2\lambda \left( d+c\right) }\text{
for }x\in \Omega _{d+c},  \label{3.190}
\end{equation}%
we obtain from (\ref{3.19})%
\begin{equation}
\left\Vert \widetilde{V}\right\Vert _{H^{1}\left( \Omega _{d+c}\right)
}^{2}\leq C_{2}e^{2\lambda m}\sigma ^{2}+C_{2}e^{-2\lambda c}\left\Vert 
\widetilde{V}\right\Vert _{H^{2}\left( \Omega \right) }^{2},\forall \lambda
\geq \lambda _{1}.  \label{3.20}
\end{equation}%
Choose $\lambda =\lambda \left( \sigma ,m,c\right) $ such that $e^{2\lambda
m}\sigma ^{2}=e^{-2\lambda c}.$ Hence, $\lambda =\ln \sigma ^{-1/\left(
m+c\right) }.$ We assume that the number $\sigma _{0}$ is so small that $\ln
\sigma _{0}^{-1/\left( m+c\right) }>\lambda _{1}.$ Hence, by (\ref{3.20})
for $\sigma \in \left( 0,\sigma _{0}\right) $ 
\begin{equation}
\left\Vert \widetilde{V}\right\Vert _{H^{1}\left( \Omega _{d+c}\right)
}^{2}\leq C_{2}\left( 1+\left\Vert \widetilde{V}\right\Vert _{H^{2}\left(
\Omega \right) }^{2}\right) \sigma ^{2\rho },\text{ }\rho =c/\left(
m+c\right) .\text{ }\square  \label{3.21}
\end{equation}

\section{Convexification}

\label{sec:4}

\subsection{Weighted Tikhonov-like functional}

\label{sec:4.1}

Assume that there exists a vector function $p\in C^{2}\left( \overline{%
\Omega }\right) $ such that 
\begin{equation}
p\mid _{\Gamma }=p_{0}\left( x\right) ,\partial _{n}p\mid _{\Gamma
}=p_{1}\left( x\right) ,  \label{4.1}
\end{equation}%
where functions $p_{0},p_{1}$ are defined in (\ref{3.11}). Consider the
vector function $W\left( x\right) ,$ 
\begin{equation}
W\left( x\right) =\left( w_{0},w_{1},...,w_{N-1}\right) ^{T}\left( x\right)
=V\left( x\right) -p\left( x\right) .  \label{4.10}
\end{equation}%
Then the problem (\ref{3.10}), (\ref{3.11}) becomes%
\begin{eqnarray}
L\left( x,p,W\right) &:&=A_{0}W-Q\left( x,\nabla p,\nabla W\right) +A_{0}p=0,
\label{4.2} \\
W &\mid &_{\Gamma }=\partial _{n}W\mid _{\Gamma }=0,  \label{4.3}
\end{eqnarray}%
where the $N-$Dimensional vector function $Q\in C^{1}\left( \mathbb{R}%
^{s_{2}}\right) ,s_{2}=n\left( 2N+1\right) $ and each component of $Q$ is a
quadratic function with respect to first derivatives $\partial
_{x_{k}}w_{i}\left( x\right) ,\partial _{x_{k}}p_{i}\left( x\right) ,$ where 
$k=1,...,n$ and $i=0,...,N-1.$

Let $s=\left[ n/2\right] +2,$ where $\left[ n/2\right] $ is the largest
integer, which does not exceed $n/2.$ Consider the space $H^{s}\left( \Omega
\right) .$ By the embedding theorem $H^{s}\left( \Omega \right) \subset
C^{1}\left( \overline{\Omega }\right) $ and with a generic constant $C>0$ 
\begin{equation}
\left\Vert f\right\Vert _{C^{1}\left( \overline{\Omega }\right) }\leq
C\left\Vert f\right\Vert _{H^{s}\left( \Omega \right) },\forall f\in
H^{s}\left( \Omega \right) .  \label{4.4}
\end{equation}%
Introduce the space $H_{0,\Gamma }^{s}\left( \Omega \right) $ of $N-$%
Dimensional vector functions $W\left( x\right) $ as%
\begin{equation*}
H_{0,\Gamma }^{s}\left( \Omega \right) =\left\{ W\in H^{s}\left( \Omega
\right) :W\mid _{\Gamma }=\partial _{n}W\mid _{\Gamma }=0\right\} .
\end{equation*}%
Let $R>0$ be an arbitrary number. Denote 
\begin{equation*}
B\left( R\right) =\left\{ W\in H_{0,\Gamma }^{s}\left( \Omega \right)
:\left\Vert W\right\Vert _{H^{s}\left( \Omega \right) }<R\right\} .
\end{equation*}%
As in Theorem 3.1, choose a number $c>0$ such that $\Omega _{d+c}\neq
\varnothing .$ Obviously $\Omega _{d+c}\subset \Omega _{d}.$ To solve the
problem (\ref{4.2}), (\ref{4.3}) numerically, consider the following
weighted Tikhonov-like functional with the CWF $\varphi _{\lambda
}^{2}\left( x\right) $ in it:%
\begin{equation}
J_{\lambda ,\gamma }\left( W\right) =e^{-2\lambda \left( d+c\right) }%
\displaystyle\int\limits_{\Omega }\left[ L\left( x,p,W\right) \right]
^{2}\varphi _{\lambda }^{2}\left( x\right) dx+\gamma \left\Vert W\right\Vert
_{H^{s}\left( \Omega \right) }^{2},  \label{4.5}
\end{equation}%
where $\gamma >0$ is the regularization parameter and the multiplier $%
e^{-2\lambda \left( d+c\right) }$ is introduced here in order to balance
first and second terms in the right hand side of (\ref{4.5}).

\textbf{Minimization Problem}. \emph{Minimize the functional }$J_{\lambda
,\gamma }\left( W\right) $\emph{\ on the closed ball }$\overline{B\left(
R\right) }.$

The second term in the right hand side of (\ref{4.5}) is taken in the norm
of the space $H^{s}\left( \Omega \right) $ in order to make sure that the
iterative terms of the gradient projection method applied to the functional $%
J_{\lambda ,\gamma }\left( W\right) $ belong to the space $C^{1}\left( 
\overline{\Omega }\right) $, see (\ref{4.4}).

\textbf{Theorem 4.1}. \emph{The functional }$J_{\lambda ,\gamma }\left(
W\right) $\emph{\ has the Frech\'{e}t derivative }$J_{\lambda ,\gamma
}^{\prime }\left( W\right) $\emph{\ at every point }$W\in H_{0,\Gamma
}^{s}\left( \Omega \right) .$\emph{\ This derivative satisfies the Lipschitz
condition in }$\overline{B\left( R\right) },$\emph{\ i.e. there exists a
constant }$Lip=Lip\left( \lambda ,\gamma ,Z,R\right) >0$\emph{\ depending
only on listed parameters such that for all }$\lambda ,\gamma >0$ 
\begin{equation*}
\left\Vert J_{\lambda ,\gamma }^{\prime }\left( W_{1}\right) -J_{\lambda
,\gamma }^{\prime }\left( W_{2}\right) \right\Vert _{H^{s}\left( \Omega
\right) }\leq Lip\left\Vert W_{1}-W_{2}\right\Vert _{H^{s}\left( \Omega
\right) },\text{ }\forall W_{1},W_{2}\in \overline{B\left( R\right) }.
\end{equation*}

\textbf{Theorem 4.2} (global strict convexity). \emph{Choose a number }$D>0$%
\emph{\ such that }$\left\Vert p\right\Vert _{C^{2}\left( \overline{\Omega }%
\right) }\leq D.$ \emph{There exists a sufficiently large number }$\lambda
_{2}=\lambda _{2}\left( \Omega ,R,Z,D,d,\xi \right) $\emph{\ such that for
all }$\lambda \geq \lambda _{2}$\emph{\ and for }$\gamma \in \left[
e^{-\lambda c},1\right) $\emph{\ the functional }$J_{\lambda ,\gamma }\left(
W\right) $\emph{\ is strictly convex on }$B\left( R\right) $\emph{, i.e. } 
\begin{eqnarray}
J_{\lambda ,\gamma }\left( W_{2}\right) -J_{\lambda ,\gamma }\left(
W_{1}\right) -J_{\lambda ,\gamma }^{\prime }\left( W_{1}\right) \left(
W_{2}-W_{1}\right) &\geq &  \label{4.6} \\
C_{1}\left\Vert W_{2}-W_{1}\right\Vert _{H^{1}\left( \Omega _{d+c}\right)
}^{2}+\frac{\gamma }{2}\left\Vert W_{2}-W_{1}\right\Vert _{H^{s}\left(
\Omega \right) }^{2},\text{ \ }\forall W_{1},W_{2} &\in &\overline{B\left(
R\right) }.  \notag
\end{eqnarray}

\textbf{Remark 4.1.} \emph{Since the regularization parameter }$\gamma \in
\left( e^{-\lambda c},1\right) $\emph{, then this allows values of }$\gamma $
\emph{to be small. Also, the presence of the first term in the right hand
side of (\ref{4.6}) indicates that the stable reconstruction should be
expected in the subdomain }$\Omega _{d+c}$\emph{\ rather than in the whole
domain }$\Omega .$ \emph{Theorem 4.4 confirms the latter.}

Let $P_{\overline{B\left( R\right) }}:H_{0,\Gamma }^{s}\left( \Omega \right)
\rightarrow \overline{B\left( R\right) }$ be the projection operator of the
Hilbert space $H_{0,\Gamma }^{s}\left( \Omega \right) $ on the closed ball $%
\overline{B\left( R\right) }.$ Let $\varsigma \in \left( 0,1\right) $ be a
number which we will choose later. Let $W_{0}\in B\left( R\right) $ be an
arbitrary point. The gradient projection method of the minimization of the
functional $J_{\lambda ,\gamma }\left( W\right) $ on the set $\overline{%
B\left( R\right) }$ is defined by the following sequence:%
\begin{equation}
W_{n}=P_{\overline{B\left( R\right) }}\left( W_{n-1}-\varsigma J_{\lambda
,\gamma }^{\prime }\left( W_{n-1}\right) \right) ;n=1,2,...  \label{4.7}
\end{equation}

\textbf{Theorem 4.3}. \emph{Let }$\lambda _{2}=\lambda _{2}\left( \Omega
,R,Z,D,c,d,\xi \right) $\emph{\ be the number introduced in Theorem 4.2. Fix
a number }$\lambda \geq $\emph{\ }$\lambda _{2}$\emph{\ and let the
regularization parameter }$\gamma \in \left[ e^{-\lambda c},1\right) .$\emph{%
\ Then there exists unique minimizer }$W_{\min }\in \overline{B\left(
R\right) }$\emph{\ of the functional }$J_{\lambda ,\gamma }\left( W\right) $%
\emph{\ on the set }$\overline{B\left( R\right) }.$\emph{\ Furthermore,
there exists a sufficiently small number }$\varsigma _{0}=\varsigma
_{0}\left( \Omega ,R,Z,D,c,d,\xi ,\lambda \right) \in \left( 0,1\right) $%
\emph{\ depending only on listed parameters such that for every }$\varsigma
\in \left( 0,\varsigma _{0}\right) $\emph{\ there exists a number }$%
q=q\left( \varsigma \right) \in \left( 0,1\right) $\emph{\ such that the
sequence (\ref{4.7}) converges to }$W_{\min },$%
\begin{equation*}
\left\Vert W_{n}-W_{\min }\right\Vert _{H^{s}\left( \Omega \right) }\leq
q^{n}\left\Vert W_{0}-W_{\min }\right\Vert _{H^{s}\left( \Omega \right)
},n=1,2,...
\end{equation*}

Consider now the question of the convergence of the sequence (\ref{4.7}) to
the exact solution $W^{\ast }$ of the problem\emph{\ }(\ref{4.2}), (\ref{4.3}%
).

\textbf{Theorem 4.4}. \emph{Assume that there exists exact solution }$%
W^{\ast }\in B\left( R\right) $\emph{\ of the problem (\ref{4.2}), (\ref{4.3}%
) with the exact data }$p^{\ast }\in C^{2}\left( \overline{\Omega }\right) .$%
\emph{\ Let }$p\in C^{2}\left( \overline{\Omega }\right) $\emph{\ be the
noisy data. Assume that }$\left\Vert p-p^{\ast }\right\Vert _{C^{2}\left( 
\overline{\Omega }\right) }\leq \sigma ,$\emph{\ where }$\sigma \in \left(
0,1\right) $\emph{\ is the level of the error in the data. Also, assume that
the }$C^{2}\left( \overline{\Omega }\right) -$\emph{norm of the exact data }$%
p^{\ast }$\emph{\ is bounded by an a priori given constant }$M^{\ast },$%
\emph{\ i.e. }$\left\Vert p^{\ast }\right\Vert _{C^{2}\left( \overline{%
\Omega }\right) }\leq M^{\ast }$ (\emph{then} $\left\Vert p\right\Vert
_{C^{2}\left( \overline{\Omega }\right) }\leq M^{\ast }+1).$\emph{\ Let }$%
\lambda _{2}=\lambda _{2}\left( \Omega ,R,Z,D,c,d,\xi \right) $\emph{\ be
the number of Theorem 4.2. Then there exists a number }$\lambda _{3}=\lambda
_{3}\left( \Omega ,R,Z,D,c,d,\xi ,M^{\ast }\right) >\lambda _{2},$\emph{\ a
sufficiently small number }$\sigma _{1}=\sigma _{1}\left( \Omega
,R,Z,D,c,d,\xi ,M^{\ast }\right) \in \left( 0,1\right) $\emph{\ and a number 
}$\theta =c/\left( 8m\right) \in \left( 0,1\right) $\emph{\ such that if }$%
\ln \sigma _{1}^{-2\theta /c}>\lambda _{3},$ \emph{then if for any\ }$\sigma
\in \left( 0,\sigma _{1}\right) $\emph{\ one chooses }$\lambda =\ln \sigma
^{-2\theta /c}$\emph{\ and }$\gamma =e^{-\lambda c}=\sigma ^{2\theta },$%
\emph{\ then the following convergence estimate holds for the sequence (\ref%
{4.7}): }%
\begin{equation}
\left\Vert W^{\ast }-W_{n}\right\Vert _{H^{1}\left( \Omega _{d+c}\right)
}\leq C_{4}\sigma ^{\theta }+q^{n}\left\Vert W_{0}-W_{\min }\right\Vert
_{H^{s}\left( \Omega \right) },n=1,2,...,  \label{4.8}
\end{equation}%
\emph{where \ the number }$q=q\left( \varsigma \right) \in \left( 0,1\right) 
$\emph{\ and }$\varsigma \in \left( 0,\varsigma _{1}\right) ,$\emph{\ where }%
$\varsigma _{1}=\varsigma _{1}\left( \Omega ,R,Z,D,c,d,\xi ,M^{\ast }\right)
\in \left( 0,1\right) $\emph{\ is a sufficiently small number. In (\ref{4.8}%
) }$C_{4}=C_{4}\left( \Omega ,R,Z,D,c,d,\xi ,M^{\ast }\right) =const.>0.$%
\emph{\ All numbers here depend only on listed parameters. }

Theorem 4.2 is the central one among theorems 4.1-4.4. Thus, we prove
Theorem 4.2 below. A similar theorem 3.2 was proven in \cite{KB} for the
case of a single quasilinear PDE, as opposed to the coupled system (\ref{4.2}%
) of quasilinear PDEs. As to the rest of theorems of this section, we omit
their proofs referring the reader to proofs of similar theorems in \cite{KB}
for the case of a single quasilinear PDE.

\textbf{Remark 4.2}. \emph{Theorem 4.2 means the convexification. Unlike (%
\ref{4.8}), in the case of non-convex functionals there is no guarantee that
a gradient-like method converges to the exact solution starting from an
arbitrary point. Since the starting point }$W_{0}\in B\left( R\right) $\emph{%
\ of the iterative process (\ref{4.7}) is an arbitrary one and since
smallness restrictions on the radius }$R$\emph{\ are not imposed, then
convergence estimate (\ref{4.8}) means the \textbf{global convergence} in
the space }$H^{1}\left( \Omega _{d+c}\right) .$

\textbf{Proof of Theorem 4.2}. In this proof, $C_{3}=C_{3}\left( \Omega
,R,Z,D,c,d,\xi \right) >0$ denotes different constants depending only on
listed parameters. Let $W_{1}$,$W_{2}\in \overline{B\left( R\right) }$ be
two arbitrary functions. Denote $W_{2}-W_{1}=h=\left( h_{0}\left( x\right)
,...,h_{N-1}\left( x\right) \right) ^{T}.$ Since each component of the
vector function $Q\left( x,\nabla p,\nabla W\right) $ in (\ref{4.2}) is a
quadratic function with respect to first derivatives $\partial
_{x_{k}}w_{i}\left( x\right) ,\partial _{x_{k}}p_{i}\left( x\right) ,$ we
have 
\begin{equation}
Q\left( x,\nabla p,\nabla W_{1}+\nabla h\right) =  \label{5.1}
\end{equation}%
\begin{equation*}
Q\left( x,\nabla p,\nabla W_{1}\right) +Q^{\left( 1\right) }\left( x,\nabla
p,\nabla W_{1},\nabla h\right) +Q^{\left( 2\right) }\left( x,\nabla p,\nabla
W_{1},\nabla h\right) .
\end{equation*}%
Here, each component of the vector function $Q^{\left( 1\right) }$ is linear
with respect to derivatives $\partial _{x_{k}}h_{i}$ and each component of
the vector function $Q^{\left( 2\right) }$ contains only quadratic terms $%
\left( \partial _{x_{k}}h_{i}\right) \cdot \left( \partial
_{x_{l}}h_{j}\right) .$ Hence, the following estimates hold for all $x\in 
\overline{\Omega }:$ 
\begin{equation}
\left\vert Q^{\left( 1\right) }\left( x,\nabla p,\nabla W_{1},\nabla
h\right) \right\vert \leq C_{3}\left\vert \nabla h\right\vert ,\text{ }%
\left\vert Q^{\left( 2\right) }\left( x,\nabla p,\nabla W_{1},\nabla
h\right) \right\vert \leq C_{3}\left\vert \nabla h\right\vert ^{2}.
\label{5.2}
\end{equation}%
By (\ref{4.2}) and (\ref{5.1}) 
\begin{equation*}
\left[ L\left( x,p,W_{1}+h\right) \right] ^{2}-\left[ L\left(
x,p,W_{1}+h\right) \right] ^{2}=Lin\left( x,p,h\right)
\end{equation*}%
\begin{equation}
+\left( A_{0}\left( h\right) \right) ^{2}+2A_{0}\left( h\right) \left[
Q^{\left( 1\right) }\left( x,\nabla p,\nabla W_{1},\nabla h\right)
+Q^{\left( 2\right) }\left( x,\nabla p,\nabla W_{1},\nabla h\right) \right]
\label{5.3}
\end{equation}%
\begin{equation*}
+\left[ Q^{\left( 1\right) }\left( x,\nabla p,\nabla W_{1},\nabla h\right)
+Q^{\left( 2\right) }\left( x,\nabla p,\nabla W_{1},\nabla h\right) \right]
^{2},
\end{equation*}%
where the functional $Lin\left( x,p,h\right) $ depends linearly on $h$.
Combining (\ref{5.3}) with the Cauchy-Schwarz inequality and as well as with
(\ref{5.2}), we obtain 
\begin{equation}
\left[ L\left( x,p,W_{1}+h\right) \right] ^{2}-\left[ L\left(
x,p,W_{1}+h\right) \right] ^{2}-Lin\left( x,p,h\right) \geq \frac{1}{2}%
\left( A_{0}\left( h\right) \right) ^{2}-C_{3}\left( \nabla h\right) ^{2}.
\label{5.4}
\end{equation}%
Let $\left\{ ,\right\} $ be the scalar product in the space of such real
valued $N-$dimensional vector functions whose components belong to $%
H^{s}\left( \Omega \right) .$ Then (\ref{4.5}) and (\ref{5.4}) imply that 
\begin{equation*}
J_{\lambda ,\gamma }\left( W_{1}+h\right) -J_{\lambda ,\gamma }\left(
W_{1}\right) -e^{-2\lambda \left( d+c\right) }\displaystyle%
\int\limits_{\Omega }Lin\left( x,p,h\right) \varphi _{\lambda
}^{2}dx-2\gamma \left\{ W,h\right\}
\end{equation*}%
\begin{equation}
\geq \frac{1}{2}\displaystyle\int\limits_{\Omega }\left( A_{0}\left(
h\right) \right) ^{2}\varphi _{\lambda }^{2}dx-C_{3}\displaystyle%
\int\limits_{\Omega }\left( \nabla h\right) ^{2}\varphi _{\lambda
}^{2}dx+\gamma \left\Vert h\right\Vert _{H^{s}\left( \Omega \right) }^{2}.
\label{5.5}
\end{equation}%
It easily follows from the proof of theorem 3.1 of \cite{KB}, which is a
close analog of Theorem 4.1, that 
\begin{equation}
J_{\lambda ,\gamma }^{\prime }\left( W_{1}\right) \left( h\right)
=e^{-2\lambda \left( d+c\right) }\displaystyle\int\limits_{\Omega }Lin\left(
x,p,h\right) \varphi _{\lambda }^{2}dx+2\gamma \left\{ W,h\right\} .
\label{5.6}
\end{equation}%
Applying the Carleman estimate (\ref{2.90})-(\ref{2.92}) to the right hand
side of (\ref{5.5}), we obtain for $\lambda \geq \lambda _{0}:$%
\begin{equation*}
\frac{e^{-2\lambda \left( d+c\right) }}{2}\displaystyle\int\limits_{\Omega
}\left( A_{0}\left( h\right) \right) ^{2}\varphi _{\lambda
}^{2}dx-C_{3}e^{-2\lambda \left( d+c\right) }\displaystyle%
\int\limits_{\Omega }\left( \nabla h\right) ^{2}\varphi _{\lambda
}^{2}dx+\gamma \left\Vert h\right\Vert _{H^{s}\left( \Omega \right) }^{2}\geq
\end{equation*}%
\begin{equation*}
\frac{e^{-2\lambda \left( d+c\right) }}{2}\displaystyle\int\limits_{\Omega
_{d}}\left( A_{0}\left( h\right) \right) ^{2}\varphi _{\lambda
}^{2}dx-C_{3}e^{-2\lambda \left( d+c\right) }\displaystyle%
\int\limits_{\Omega _{d}}\left( \nabla h\right) ^{2}\varphi _{\lambda
}^{2}dx-C_{3}\displaystyle\int\limits_{\Omega \diagdown \Omega _{d}}\left(
\nabla h\right) ^{2}\varphi _{\lambda }^{2}dx+\gamma \left\Vert h\right\Vert
_{H^{s}\left( \Omega \right) }^{2}
\end{equation*}%
\begin{equation*}
\geq C_{1}e^{-2\lambda \left( d+c\right) }\lambda \displaystyle%
\int\limits_{\Omega _{d}}\left( \nabla h\right) ^{2}\varphi _{\lambda
}^{2}dx+C_{1}e^{-2\lambda \left( d+c\right) }\lambda ^{3}\displaystyle%
\int\limits_{\Omega _{d}}h^{2}\varphi _{\lambda }^{2}dx-C_{3}e^{-2\lambda
\left( d+c\right) }\displaystyle\int\limits_{\Omega _{d}}\left( \nabla
h\right) ^{2}\varphi _{\lambda }^{2}dx
\end{equation*}%
\begin{equation}
-C_{3}e^{-2\lambda \left( d+c\right) }\displaystyle\int\limits_{\Omega
\diagdown \Omega _{d}}\left( \nabla h\right) ^{2}\varphi _{\lambda
}^{2}dx-C_{1}\lambda ^{3}e^{-2\lambda c}\displaystyle\int\limits_{\xi
_{d}}\left( \left( \nabla h\right) ^{2}+h^{2}\right) dS+\gamma \left\Vert
h\right\Vert _{H^{s}\left( \Omega \right) }^{2}.  \label{5.7}
\end{equation}

Choose $\lambda _{2}=\lambda _{2}\left( \Omega ,R,Z,D,d,\xi \right) \geq
\lambda _{0}$ so large that $C_{1}\lambda _{2}/2>C_{3}.$ Also, observe that 
\begin{equation*}
\varphi _{\lambda }^{2}\left( x\right) \leq e^{2\lambda d},\forall x\in
\Omega \diagdown \Omega _{d}\text{ and }\left\Vert \nabla h\right\Vert
_{L_{2}\left( \xi _{d}\right) }^{2}+\left\Vert h\right\Vert _{L_{2}\left(
\xi _{d}\right) }^{2}\leq C_{3}\left\Vert h\right\Vert _{H^{s}\left( \Omega
\right) }^{2}.
\end{equation*}%
Hence, taking into account (\ref{3.190}), we obtain from (\ref{5.7})%
\begin{equation}
\frac{e^{-2\lambda \left( d+c\right) }}{2}\displaystyle\int\limits_{\Omega
}\left( A_{0}\left( h\right) \right) ^{2}\varphi _{\lambda
}^{2}dx-C_{3}e^{-2\lambda \left( d+c\right) }\displaystyle%
\int\limits_{\Omega }\left( \nabla h\right) ^{2}\varphi _{\lambda
}^{2}dx+\gamma \left\Vert h\right\Vert _{H^{s}\left( \Omega \right) }^{2}
\label{5.8}
\end{equation}%
\begin{equation*}
\geq C_{1}\left\Vert h\right\Vert _{H^{1}\left( \Omega _{d+c}\right)
}^{2}+\left( \gamma -C_{3}e^{-2\lambda c}\right) \left\Vert h\right\Vert
_{H^{s}\left( \Omega \right) }^{2},\forall \lambda \geq \lambda _{2}.
\end{equation*}%
Since $\gamma \in \left[ e^{-\lambda c},1\right) ,$ then (\ref{5.5}), (\ref%
{5.6}) and (\ref{5.8}) imply (\ref{4.6}). $\square $

\subsection{Numerical scheme}

\label{sec:4.2}

The numerical scheme for the above technique is as follows:

\textbf{Step 1}. Using symbolic computations for the Gram-Schmidt
Orthonormalization procedure in $L_{2}\left( 0,1\right) $, obtain functions $%
\left\{ \psi _{m}\left( x_{0}\right) \right\} _{m=0}^{N-1},x_{0}\in \left(
0,1\right) $ from functions $\left\{ x_{0}^{m}e^{x_{0}}\right\} _{m=0}^{N-1}$
for a reasonable integer $N\geq 2.$ Alternatively, if $x_{0}\in \left(
0,\infty \right) ,$ then use Laguerre functions.

\textbf{Step 2}. Sequentially obtain problems (\ref{3.3}), (\ref{3.4}), then
(\ref{3.10}), (\ref{3.11}) and then (\ref{4.2}), (\ref{4.3}) for the
specific operator $A$.

\textbf{Step 3}. Minimize the functional (\ref{4.5}) on the set $\overline{%
B\left( R\right) }$ using the gradient projection method.

\textbf{Step 4}.\ Let $W_{\min }\left( x\right) $ be the minimizer of the
functional (\ref{4.5}) on the set $\overline{B\left( R\right) }$ (Theorem
4.3). Set $V_{\min }\left( x\right) =W_{\min }\left( x\right) +p\left(
x\right) .$ Next, use the first formula (\ref{3.8}), then use (\ref{3.5})
and finally use (\ref{3.1}).

\section{Elliptic Equation}

\label{sec:5}

The goal of sections 5-8 is to provide some specific examples of CIPs for
which the above technique works. We point out that a \emph{variety} of other
examples are possible.

\subsection{The general case}

\label{sec:5.1}

In this case conditions (\ref{2.1}), (\ref{2.2}) are specified as%
\begin{eqnarray}
Au &=&\displaystyle\sum \limits_{i,j=1}^{n}a_{i,j}\left( x\right)
u_{x_{i}x_{j}}+\displaystyle\sum \limits_{j=1}^{n}b_{j}\left( x\right)
u_{x_{j}}+a_{0}\left( x\right) u,x\in \mathbb{R}^{n},  \label{6.1} \\
A_{0}u &=&\displaystyle\sum \limits_{i,j=1}^{n}a_{i,j}\left( x\right)
u_{x_{i}x_{j}},  \label{6.2}
\end{eqnarray}%
where $a_{i,j}\left( x\right) =a_{j,i}\left( x\right) ,\forall i,j.$ We
assume that obvious analogs of conditions (\ref{2.30})-(\ref{2.4}) are
valid. Also, we assume that there exist two constants $\mu _{1},\mu
_{2}>0,\mu _{1}\leq \mu _{2}$ such that%
\begin{equation}
\mu _{1}\left\vert \eta \right\vert ^{2}\leq \displaystyle\sum
\limits_{i,j=1}^{n}a_{i,j}\left( x\right) \eta _{i}\eta _{j}\leq \mu
_{2}\left\vert \eta \right\vert ^{2},\forall x\in \mathbb{R}^{n},\forall
\eta =\left( \eta _{1},...\eta _{n}\right) \in \mathbb{R}^{n}.  \label{6.20}
\end{equation}%
Let the domain $\Omega \subset \left\{ x_{n}>0\right\} $ and $\overline{x}%
^{0}=\left( 0,...,0,-1\right) .$\ Let $I$ be the interval defined in (\ref%
{2.10}). Hence, (\ref{2.11}) is valid. Choose a number $\omega >1$ and
assume that 
\begin{equation}
\Gamma =\left\{ x_{n}=0,\left( x_{1}-1/2\right) ^{2}/\omega ^{2}+%
\displaystyle\sum \limits_{k=2}^{n-1}x_{k}^{2}<\frac{1}{4}\right\} \subset
\partial \Omega .  \label{6.3}
\end{equation}%
Hence, $\left\{ x_{1}\in \left( 0,1\right) ,\overline{x}=0\right\} \subset
\left( \Gamma \cap \left\{ x_{2}=...=x_{n-1}=0\right\} \right) ,$ which
means that the interval $I$, over which the point source is running, is
\textquotedblleft observable" from the hypersurface $\Gamma ,$ where our
Dirichlet and Neumann data are given, see Condition 4 in section 2.2. This
\textquotedblleft observability" seems to be important for successful
numerical studies.

Suppose that a distribution $u\left( x,x_{0}\right) $ satisfies equation (%
\ref{2.12}) and Conditions 1-4 of section 2.2. Define the functions $\xi
\left( x\right) $ and the CWF $\varphi _{\lambda }\left( x\right) $ as%
\begin{equation}
\xi \left( x\right) =\left[ x_{n}+\left( x_{1}-1/2\right) ^{2}/\omega ^{2}+%
\displaystyle\sum \limits_{k=2}^{n-1}x_{k}^{2}+\frac{1}{4}\right] ^{-\nu },%
\text{ }\varphi _{\lambda }\left( x\right) =\exp \left( \lambda \xi \left(
x\right) \right) ,  \label{6.4}
\end{equation}%
where $\nu =\nu \left( \omega \right) >1$ is a parameter depending on $%
\omega .$ Assume that $\left\{ \xi \left( x\right) >2^{\nu },x_{n}>0\right\}
=\Omega _{2^{\nu }}\subset \Omega .$ Then by (\ref{2.7}) (\ref{6.3}) and (%
\ref{6.4}) $\Gamma =\Gamma _{2^{\nu }}.$ It follows from results of \S 1 of
Chapter 4 of \cite{LRS} that a direct analog of the Carleman estimate of (%
\ref{2.90})-(\ref{2.92}) is valid for the operator $A_{0}$ in (\ref{6.2})
with conditions (\ref{2.31}) and (\ref{6.20}) and with the CWF (\ref{6.4}),
as long as the parameter $\nu $ is sufficiently large.

Therefore, the above construction works in this case. The unknown
coefficient $a_{0}\left( x\right) $ can be H\"{o}lder-stable reconstructed
numerically by the above method in the domain $\Omega _{2^{\nu }+c},c>0,$ as
long as $\Omega _{2^{\nu }+c}\neq \varnothing .$ On the other hand,
uniqueness is guaranteed in the entire domain $\Omega .$ This can be proven
similarly with the conventional proof of the uniqueness of the problem of
the continuation of solutions of elliptic equations.

\subsection{Helmholtz equation}

\label{sec:5.2}

We now specify the discussion of section 6.1 for the case of the Helmholtz
equation, since it is interesting for many applications. Let $x\in \mathbb{R}%
^{3}.$ Let the function $c\left( x\right) $ be smooth in $\mathbb{R}^{3}$
and be such that $c\left( x\right) \geq 1$ in $\mathbb{R}^{3}$ and $c\left(
x\right) =1$ for $x\in \mathbb{R}^{3}\diagdown \Omega .$ The Helmholtz
equation for the function $u\left( x,k\right) $ with the radiation condition
is and with the source located at $\left( x_{0},\overline{x}^{0}\right)
=\left( x_{0},0,-1\right) =y\left( x_{0}\right) $ is%
\begin{eqnarray}
\Delta u+k^{2}c\left( x\right) u &=&-\delta \left( x_{1}-x_{0}\right) \delta
\left( x_{2},x_{3}+1\right) ,  \label{6.5} \\
\lim_{r\rightarrow \infty }\left[ r\left( \partial _{r}u-iku\right) \right]
&=&0,r=\left\vert x\right\vert ,  \label{6.6}
\end{eqnarray}%
where $k$ is the wavenumber. In applications to scattering of
electromagnetic waves, $c\left( x\right) $ is a spatially distributed
dielectric constant. Using a comparison with the solution of the Maxwell's
equations, it was demonstrated numerically in \cite{BMM} that a simplified
model of propagation of electric wave field, based on (\ref{6.5}), (\ref{6.6}%
), can be used instead of the Maxwell's equations. This conclusion was
confirmed by many accurate imaging results of the author with coauthors,
using experimentally measured microwave data, see \cite{Liem,Liem2} for the
frequency domain data and references cited there for the time domain data.

Let $a$ and $B$ be two numbers, where $a>1,B>0$. We set 
\begin{equation}
\Omega =\left\{ x:-a<x_{1},x_{2}<a,x_{3}\in \left( 0,B\right) \right\} .
\label{6.60}
\end{equation}%
We consider two Coefficient Inverse Scattering Problems (CISPs): one with
backscattering data and another one with transmitted data.

\textbf{Coefficient Inverse Scattering Problem 1} (CISP1, backscattering
data). \emph{Let }$\Gamma _{b}=\left\{ x:-a<x_{1},x_{2}<a,x_{3}=0\right\} .$%
\emph{\ Determine the unknown coefficient }$c\left( x\right) $\emph{\ for }$%
x\in \Omega ,$\emph{\ assuming that the following two functions }$%
g_{0}\left( x_{1},x_{2},x_{0}\right) ,g_{1}\left( x_{1},x_{2},x_{0}\right) $%
\emph{\ are known for a fixed value of the wavenumber }$k=k_{0}$\emph{:}%
\begin{equation}
u\left( x,k_{0}\right) \mid _{\Gamma _{b}=}g_{0}\left(
x_{1},x_{2},x_{0}\right) ,\partial _{n}u\left( x,k_{0}\right) \mid _{\Gamma
_{b}=}g_{1}\left( x_{1},x_{2},x_{0}\right) ,\forall x_{0}\in \left[ 0,1%
\right] .  \label{6.7}
\end{equation}

\textbf{Coefficient Inverse Scattering Problem 2} (CISP2, transmitted data). 
\emph{Let }$\Gamma _{tr}=\left\{ x:-a<x_{1},x_{2}<a,x_{3}=B\right\} .$\emph{%
\ Determine the unknown coefficient }$c\left( x\right) $\emph{\ for }$x\in
\Omega ,$\emph{\ assuming that the following two functions }$g_{0}\left(
x_{1},x_{2},x_{0}\right) ,g_{1}\left( x_{1},x_{2},x_{0}\right) $\emph{\ are
known for a fixed value of the wavenumber }$k=k_{0}$\emph{:}%
\begin{equation}
u\left( x,k_{0}\right) \mid _{\Gamma _{tr}=}g_{0}\left(
x_{1},x_{2},x_{0}\right) ,\partial _{n}u\left( x,k_{0}\right) \mid _{\Gamma
_{tr}=}g_{1}\left( x_{1},x_{2},x_{0}\right) ,\forall x_{0}\in \left[ 0,1%
\right] .  \label{6.8}
\end{equation}

CISP1 has direct applications in imaging of land mines and improvised
explosive devices (IEDs) \cite{Liem,Liem2}. As to the CISP2, it has direct
applications in crosswell imaging, see Figures 1a),b).

\begin{figure}[h]
\begin{center}
\subfloat[]{\includegraphics[width=0.5\textwidth]{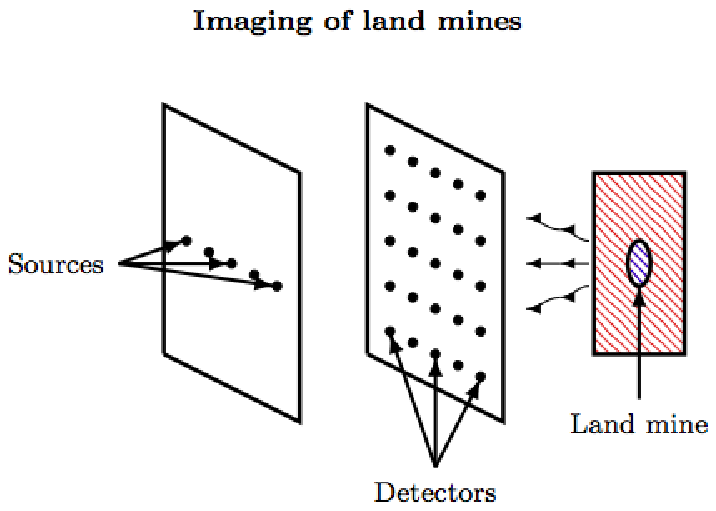}} \subfloat[]{%
\includegraphics[width=0.5\textwidth]{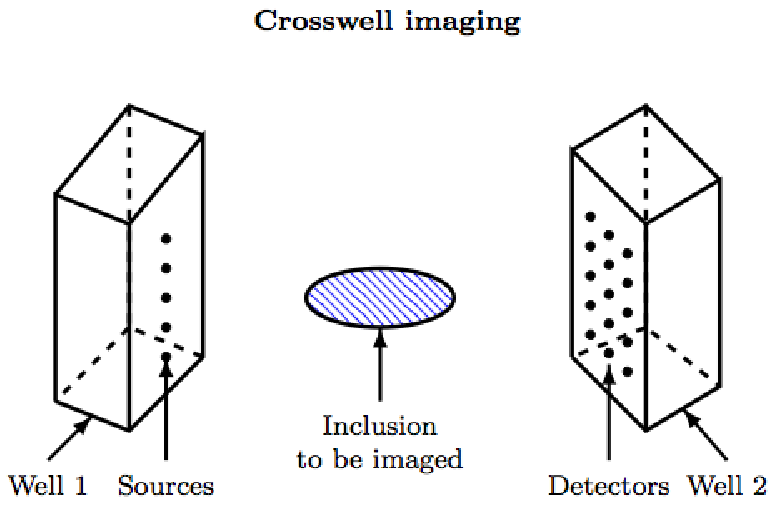}}
\end{center}
\caption{Schematic measurement diagrams for two applications. a) Imaging of
land mines and improvised explosive devices (CISP1). b) Crosswell imaging
(CISP2).}
\label{fig 1}
\end{figure}

As to the first application, the author has many publications on this
subject, in which values of dielectric constants and locations of objects
mimicking land mines and IEDs were accurately imaged, including many cases
of real data, see, e.g. \cite{Liem,Liem2} and references cited there. In
these references a globally convergent numerical method was applied.
However, a single location of the source and an interval of wavenumbers were
used in these publications, unlike the current case of the restricted DN
data and a single wavenumber.

\textbf{Remark 5.1}. \emph{The question \textquotedblleft How do we know
both functions }$g_{0}$\emph{\ and }$g_{1}$\emph{\ in (\ref{6.7}), (\ref{6.8}%
) if only the function }$u$ \emph{is} \emph{usually actually measured in
experiments?" is addressed in \cite{Liem,Liem2} via the so-called
\textquotedblleft data propagation procedure", see a detailed description in 
\cite{Liem2}. In the case of backscattering data, this
procedure\textquotedblleft moves" the data from the plane }$\left\{
x_{3}=-z=const.<0\right\} ,$\emph{\ where the data are originally collected,
to the surface }$\Gamma _{b}\subset \left\{ x_{3}=0\right\} ,$\emph{\ which
is closer to the targets to be imaged. A similar procedure can be arranged
for transmitted data in the case of CISP2.}

It was established in \cite{KR1} that, under certain conditions imposed on
the coefficient $c\left( x\right) ,$ the following asymptotic expansion of
the function $u\left( x,x_{0},k\right) $ takes place%
\begin{equation}
u\left( x,x_{0},k\right) =A\left( x,y\right) e^{ik\tau \left( x,y\right)
}\left( 1+h\left( x,y\left( x_{0}\right) ,k\right) \right) ,k\rightarrow
\infty ,\forall x\in \overline{\Omega },\forall x_{0}\in \left[ 0,1\right] ,
\label{6.80}
\end{equation}%
where the function $h\left( x,y\left( x_{0}\right) ,k\right) $ is such that $%
\left\vert h\left( x,y\left( x_{0}\right) ,k\right) \right\vert \leq
C_{5}/k,\forall x\in \overline{\Omega },$ where $C_{5}>0$ is a certain
constant independent on $x,k$. In (\ref{6.80}) $A\left( x,y\left(
x_{0}\right) \right) >0,\forall x\in \overline{\Omega },\forall x_{0}\in %
\left[ 0,1\right] $ and $\tau \left( x,y\left( x_{0}\right) \right) $ is the
travel time between points $y\left( x_{0}\right) $ and $x$ due to the
Riemannian metric generated by the function $c\left( x\right) .$ It follows
from (\ref{6.80}) that $u\left( x,x_{0},k\right) \neq 0,\forall x\in 
\overline{\Omega },\forall x_{0}\in \left[ 0,1\right] $ for sufficiently
large values of $k$. This replaces Condition 3 of section 2.2 for the case
of Helmholtz equation. Furthermore, $\log u\left( x,x_{0},k\right) $ can be
uniquely defined as:%
\begin{equation*}
\log u\left( x,x_{0},k\right) =ik\tau \left( x,y\left( x_{0}\right) \right)
+\ln A\left( x,y\left( x_{0}\right) \right) +\log \left( 1+h\left( x,y\left(
x_{0}\right) ,k\right) \right)
\end{equation*}%
\begin{equation*}
:=ik\tau \left( x,y\left( x_{0}\right) \right) +\ln A\left( x,y\left(
x_{0}\right) \right) +\displaystyle\sum\limits_{n=1}^{\infty }\frac{\left(
-1\right) ^{n-1}}{n}h^{n}\left( x,y\left( x_{0}\right) ,k\right) ,\text{ }%
x\in \overline{\Omega }.
\end{equation*}%
In the case of CISP1 the CWF (\ref{6.4}) can be chosen. In the case of CISP2
an obvious modification of (\ref{6.4}) can be chosen. We conclude,
therefore, that the above technique is applicable to both CISP1 and CISP2,
as long as $k_{0}$ is sufficiently large.

\textbf{Remark 5.2}. \emph{Note that even though the theory of the method of 
\cite{Liem,Liem2} also works only for sufficiently large values of }$k$\emph{%
, successful numerical studies for real data were conducted in \cite%
{Liem,Liem2} for quite reasonable values of }$k$\emph{: }$k\in \left[
6.25,6.7\right] $ \emph{in} \emph{\cite{Liem} and }$k\in \left[ 5.31,5.69%
\right] $ \emph{in} \emph{\cite{Liem2}, see \cite{Liem,Liem2} for
explanations of these choices from the standpoint of Physics. This indicates
that the technique of the current paper can actually work for CISP1 and
CISP2 for reasonable values of }$k_{0}.$

\subsection{Electrical impedance tomography (EIT)}

\label{sec:5.3}

In this case we consider CIP2 of section 2.2. Below for any $\alpha \in
\left( 0,1\right) $ and any integer $k\geq 0$ let $C^{k+\alpha }$ and $%
C^{2k+\alpha ,k+\alpha /2}$ be the H\"{o}lder spaces for elliptic and
parabolic equations respectively. Let the function $\sigma \in C^{2+\alpha
}\left( \overline{G}\right) ,\sigma \left( x\right) \geq \sigma
_{0}=const.>0,\forall x\in G$ and $\sigma \left( x\right) =1$ in $G\diagdown
\Omega .$ In addition, let $\partial G\in C^{2+\alpha }.$ The boundary value
problem for EIT is%
\begin{eqnarray*}
\func{div}\left( \sigma \left( x\right) \nabla u\right) &=&-f\left(
x_{1}-x_{0}\right) \chi \left( \overline{x}-\overline{x}^{0}\right) ,\forall
x_{0}\in \left[ 0,1\right] , \\
u &\mid &_{x\in \partial G}=0,\forall x_{0}\in \left[ 0,1\right] .
\end{eqnarray*}%
Introducing the well known change of variables $v=\sqrt{\sigma }u,$ we obtain%
\begin{eqnarray}
\Delta v+a_{0}\left( x\right) v &=&-f\left( x_{1}-x_{0}\right) \chi \left( 
\overline{x}-\overline{x}^{0}\right) ,\forall x_{0}\in \left[ 0,1\right] ,
\label{6.9} \\
v &\mid &_{x\in \partial G}=0,\forall x_{0}\in \left[ 0,1\right] ,
\label{6.10}
\end{eqnarray}%
where $a_{0}\left( x\right) =\Delta \left( \sqrt{\sigma \left( x\right) }%
\right) /\sqrt{\sigma \left( x\right) }.$ Recall that $f\left( 0\right) \chi
\left( 0\right) \neq 0.$ Assume that $f\left( z\right) ,\chi \left( 
\overline{x}\right) \geq 0,\forall z\in \mathbb{R},\forall \overline{x}\in 
\mathbb{R}^{n-1}$ and that $a_{0}\left( x\right) \leq 0$ in $\Omega .$ Then
theorem 6.14 of \cite{G} guarantees that there exists unique solution $v\in
C^{2+\alpha }\left( \overline{G}\right) $ of the problem (\ref{6.9}), (\ref%
{6.10}), for every $x_{0}\in \left[ 0,1\right] $. Next, the strong maximum
principle of theorem 3.5 of \cite{G} ensures that the function $v\left(
x,x_{0}\right) $ satisfies Condition 3 of section 2.2. Thus, it follows from
results of section 5.1 that the above technique is applicable to EIT.

\section{Parabolic Equation}

\label{sec:6}

Let $T>1$ be an arbitrary number. Denote $D_{T}^{n+1}=\mathbb{R}^{n}\times
\left( 0,T\right) .$ Consider the parabolic operator in $D_{T}^{n+1}$ 
\begin{eqnarray}
Au &=&u_{t}-\displaystyle\sum \limits_{i,j=1}^{n}a_{i,j}\left( x,t\right)
u_{x_{i}x_{j}}-\displaystyle\sum \limits_{j=1}^{n}b_{j}\left( x,t\right)
u_{x_{j}}+a_{0}\left( x,t\right) u,  \label{7.1} \\
A_{0}u &=&u_{t}-\displaystyle\sum \limits_{i,j=1}^{n}a_{i,j}\left(
x,t\right) u_{x_{i}x_{j}},  \label{7.2}
\end{eqnarray}%
where $a_{0}\left( x,t\right) $ is the unknown coefficient and $%
a_{i,j}\left( x,t\right) =a_{j,i}\left( x,t\right) ,\forall i,j.$ We assume
that all coefficients of the operator (\ref{7.1}) belong to $C^{1}\left( 
\overline{D_{T}^{n+1}}\right) $ and also that the obvious analog of the
ellipticity condition (\ref{6.20}) holds.

We take the same $\overline{x}^{0}=\left( 0,...,0,-1\right) $ as in section
5.1 and let the domain $\Omega _{1}\subset \left\{ x_{n}>0\right\} $. We
assume that $Au=u_{t}-\Delta u$ for $x\in \mathbb{R}^{n}\diagdown \Omega
_{1}.$ In the case of the parabolic equation the point source runs over $%
I_{1}=I\times \left\{ t=0\right\} ,$ where the interval $I$ is defined in (%
\ref{2.10}). Consider the fundamental solution $u\left( x,t,x_{0}\right) $
of the operator $A$, 
\begin{equation}
Au=\delta \left( x_{1}-x_{0},\overline{x}-\overline{x}^{0},t\right) ,\left(
x,t\right) \in D_{T}^{n+1},  \label{7.3}
\end{equation}%
\begin{equation}
u\left( x,0,x_{0}\right) =0,\forall x_{0}\in \left[ 0,1\right] .  \label{7.4}
\end{equation}%
It is well known that there exists unique solution $u\left( x,t,x_{0}\right)
\in C^{2+\alpha ,1+\alpha /2}\left( \overline{D_{T}^{n+1}}\diagdown
I_{1\varepsilon }\right) ,\forall \varepsilon >0,\forall x_{0}\in \left[ 0,1%
\right] $ of the problem (\ref{7.3}), (\ref{7.4}), see chapter 4 of \cite%
{Lad}\emph{.} Here $I_{1\varepsilon }\subset \overline{D_{T}^{n+1}}$ is
defined similarly with $I_{\varepsilon }$ in section 2.2. Furthermore, $%
u\left( x,t,x_{0}\right) >0$ for $t>0,$ see theorem 11 in chapter 2 of \cite%
{Fr}\emph{. }Hence (see Condition 3 in section 2.2), we\emph{\ }define $%
\Omega =\Omega _{1}\times \left( \zeta ,T\right) ,$ where $\zeta \in \left(
0,\left( T-1\right) /2\right) $ is a sufficiently small number.

We define the hypersurface $\Gamma $ and the function $\xi \left( x,t\right) 
$ similarly with (\ref{6.3}) and (\ref{6.4}). Choose a number $\omega >1$
and set%
\begin{equation*}
\Gamma _{0}=\left\{ x\in \mathbb{R}^{n}:x_{n}=0,\left( x_{1}-1/2\right)
^{2}/\omega ^{2}+\displaystyle\sum\limits\limits_{k=2}^{n-1}x_{k}^{2}<1/4%
\right\} ,\Gamma =\Gamma _{0}\times \left( \zeta ,T\right) .
\end{equation*}%
\begin{equation}
\xi \left( x,t\right) =\left[ x_{n}+\left( x_{1}-1/2\right) ^{2}/\omega ^{2}+%
\displaystyle\sum\limits\limits_{k=2}^{n-1}x_{k}^{2}+\left( t-T/2\right)
^{2}+1/4\right] ^{-\nu }.  \label{7.6}
\end{equation}%
and $\varphi _{\lambda }\left( x,t\right) =\exp \left( \lambda \xi \left(
x,t\right) \right) $, where $\nu =\nu \left( \omega \right) >1$ is a
parameter depending on $\omega .$ We assume that $\Gamma _{0}\subset
\partial \Omega _{1}.$ Hence, $\Gamma \subset \partial \Omega .$ In
addition, we assume that $\left\{ \xi \left( x,t\right) >2^{\nu
},x_{n}>0\right\} =\Omega _{2^{\nu }}\subset \Omega .$ Note that 
\begin{equation*}
\Gamma _{2^{\nu }}=\left\{ x_{n}=0,\left( x_{1}-1/2\right) ^{2}/\omega ^{2}+%
\displaystyle\sum\limits\limits_{k=2}^{n-1}x_{k}^{2}+\left( t-T/2\right)
^{2}<1/4\right\} \subset \Gamma .
\end{equation*}
Let the restricted DN\ data be given on $\Gamma ,$%
\begin{equation}
u\left( x,t,x_{0}\right) =g_{0}\left( x,t,x_{0}\right) ,\partial _{n}u\left(
x,t,x_{0}\right) =g_{1}\left( x,t,x_{0}\right) ,\forall \left(
x,t,x_{0}\right) \in \Gamma \times \left[ 0,1\right] .  \label{7.7}
\end{equation}

Similarly with section 5.1, results of \S 1 of Chapter 4 of \cite{LRS} imply
that a direct analog of the Carleman estimate of (\ref{2.90})-(\ref{2.92})
is valid for the operator $A_{0}$ in (\ref{7.2}), as long as the parameter $%
\nu $ is sufficiently large. Hence, the above construction works in this
case. The unknown coefficient $a_{0}\left( x,t\right) $ can be H\"{o}%
lder-stable reconstructed numerically by the above method in the domain $%
\Omega _{2^{\nu }+c}$ for any $c>0$ such that $\Omega _{2^{\nu }+c}\neq
\varnothing .$ As to the uniqueness, the entire domain $\Omega _{1}\times
\left( 0,T\right) $ works: similarly with the end of section 5.1.

\section{Hyperbolic equation}

\label{sec:7}

In this case we consider CIP2 of section 2.2. Let the domains $\Omega
_{k}\subset \mathbb{R}^{3}$ be defined as $\Omega _{k}=\left\{ \left\vert
x\right\vert <k\right\} ,k=1,2,3.$ Let $I_{\varepsilon }$ be the set defined
in section 2.2. We assume that 
\begin{equation}
I_{\varepsilon }\subset \left( \Omega _{3}\diagdown \Omega _{2}\right) .
\label{8.0}
\end{equation}%
Let the function $a\left( x,t\right) \in C\left( \overline{D_{T}^{4}}\right) 
$ be such that%
\begin{equation}
a\left( x,t\right) \geq 0\text{ in }D_{T}^{4}\text{, }a\left( x,t\right) =0%
\text{ for }x\in \mathbb{R}^{3}\diagdown \Omega _{2}.  \label{8.1}
\end{equation}%
Let $f\left( z\right) ,z\in \mathbb{R}$ and $\chi \left( \overline{x}\right)
,\overline{x}\in \mathbb{R}^{n-1}$ be functions defined in section 2.2.
Recall that $f\left( 0\right) \chi \left( 0\right) \neq 0.$ We assume that 
\begin{equation}
f\left( z\right) \geq 0,\forall z\in \mathbb{R}\text{ and }\chi \left( 
\overline{x}\right) \geq 0,\forall \overline{x}\in \mathbb{R}^{n-1}.
\label{8.2}
\end{equation}

Consider the following Cauchy problem for the function $u\left(
x,t,x_{0}\right) $%
\begin{eqnarray}
u_{tt} &=&\Delta u+a\left( x,t\right) u+f\left( x_{1}-x_{0}\right) \chi
\left( \overline{x}-\overline{x}^{0}\right) ,\left( x,t\right) \in D_{T}^{4},
\label{8.3} \\
u\left( x,0,x_{0}\right) &=&u_{t}\left( x,0,x_{0}\right) =0,  \label{8.4}
\end{eqnarray}%
where $x_{0}\in \left[ 0,1\right] $ is a parameter. Then the problem (\ref%
{8.3}), (\ref{8.4}) is equivalent with: 
\begin{equation}
u\left( x,t,x_{0}\right) =\displaystyle\int\limits_{\left\vert x-\eta
\right\vert <t}\frac{f\left( \eta _{1}-x_{0}\right) \chi \left( \overline{%
\eta }-\overline{x}^{0}\right) }{4\pi \left\vert x-\eta \right\vert }d\eta +%
\displaystyle\int\limits_{\left\vert x-\eta \right\vert <t}\frac{\left(
au\right) \left( \eta ,t-\left\vert x-\eta \right\vert \right) }{4\pi
\left\vert x-\eta \right\vert }d\eta .  \label{8.5}
\end{equation}%
On can prove (see \cite{LRS} for a similar result) that the integral
equation (\ref{8.5}) can be rewritten as Volterra integral equation, whose
solution can be represented as a series, which converges absolutely and
uniformly in any subdomain $\left( G\times \left( 0,T\right) \right) \subset 
$ $D_{T}^{4}$ and for any $x_{0}\in \left[ 0,1\right] ,$ where $G\subset 
\mathbb{R}^{3}$ is an arbitrary bounded domain. This series is%
\begin{eqnarray}
u &=&\displaystyle\sum \limits_{n=0}^{\infty }u_{n},u_{0}=\displaystyle%
\int\limits_{\left\vert x-\eta \right\vert <t}\frac{f\left( \eta
_{1}-x_{0}\right) \chi \left( \overline{\eta }-\overline{x}^{0}\right) }{%
4\pi \left\vert x-\eta \right\vert }d\eta ,  \label{8.6} \\
u_{n} &=&\displaystyle\int\limits_{\left\vert x-\eta \right\vert <t}\frac{%
\left( au_{n-1}\right) \left( \eta ,t-\left\vert x-\eta \right\vert \right) 
}{4\pi \left\vert x-\eta \right\vert }d\eta ,n\geq 1.  \label{8.60}
\end{eqnarray}%
We now prove that 
\begin{equation}
u\left( x,t,x_{0}\right) \geq C_{5}T,\forall x\in \Omega _{1},\forall
x_{0}\in \left[ 0,1\right] ,\forall t\in \left( T/4,T\right) ,\forall T>20,
\label{8.7}
\end{equation}%
where the constant $C_{5}=C_{5}\left( I_{\varepsilon },f,\chi \right) >0$
depends only on listed parameters and is independent on $T$.

Indeed, let $x\in \Omega _{1}$ and $\eta \in \Omega _{3}$ be two arbitrary
points. Let $t\in \left( T/4,T\right) $ and $T>20.$ Then 
\begin{equation}
\left\vert x-\eta \right\vert \leq \left\vert x\right\vert +\left\vert \eta
\right\vert <4<t.  \label{8.8}
\end{equation}%
Since by (\ref{8.2}) $f\left( 0\right) \chi \left( 0\right) >0,$ then (\ref%
{8.7}) follows from (\ref{8.0})-(\ref{8.2}), (\ref{8.6}), (\ref{8.60}) and (%
\ref{8.8}).

We now set $\Omega =\Omega _{1}\times \left( T/4,T\right) ,$ where $T>20$.
Next, let%
\begin{equation}
\xi \left( x,t\right) =\left\vert x\right\vert ^{2}-\varrho ^{2}\left(
t-T/2\right) ^{2},\varphi _{\lambda }\left( x,t\right) =\exp \left( \lambda
\xi \left( x,t\right) \right) .  \label{8.10}
\end{equation}%
Choose any number $d\in \left( 0,1\right) .$ Next, choose $\varrho \in
\left( 4\sqrt{1-d}/T,1\right) .$ Let $\Omega _{d}=\left\{ \left( x,t\right)
:x\in \Omega _{1},\xi \left( x,t\right) >d\right\} .$ Then%
\begin{equation}
\Omega _{d}\subset \Omega ,\Omega _{d}\cap \left\{ t=T/4\right\} =\Omega
_{d}\cap \left\{ t=T\right\} =\varnothing .  \label{8.11}
\end{equation}%
Hence, we define $\Gamma $ and $\Gamma _{d}$ as 
\begin{equation}
\Gamma =\left\{ \left( x,t\right) :\left\vert x\right\vert =1,t\in \left(
T/4,T\right) \right\} ,\Gamma _{d}=\left\{ \left( x,t\right) :\left\vert
x\right\vert =1,\xi \left( x,t\right) >d\right\} .  \label{8.12}
\end{equation}%
It follows from (\ref{8.10})-(\ref{8.12}) that $\Gamma _{d}\subset \Gamma .$

Similarly with (\ref{2.130}) we define the CIP in this case as the problem
of determining the unknown coefficient $a\left( x,t\right) \in C\left( 
\overline{D_{T}^{4}}\right) $ satisfying conditions (\ref{8.1}), given
functions $g_{0}\left( x,t,x_{0}\right) ,g_{1}\left( x,t,x_{0}\right) $,
where%
\begin{equation*}
u\left( x,t,x_{0}\right) =g_{0}\left( x,t,x_{0}\right) ,\partial _{n}u\left(
x,t,x_{0}\right) =g_{1}\left( x,t,x_{0}\right) ,\forall \left(
x,t,x_{0}\right) \in \Gamma \times \left[ 0,1\right] .
\end{equation*}%
An analog of the Carleman estimate of (\ref{2.90})-(\ref{2.92}) works for
the operator $\partial _{t}^{2}-\Delta $ with the CWF $\varphi _{\lambda
}\left( x,t\right) $ given in (\ref{8.10}), see theorem 2.2.5 in \cite{KT}.
Therefore, the above construction works for this CIP. The function $a\left(
x,t\right) $ can be reconstructed numerically by the above method in $\Omega
_{d+c}$ for any $c\in \left( 0,1-d\right) .$

\section{Some Numerical Considerations}

\label{sec:8}

We discuss in this section some practical ideas for the numerical
implementation of the procedure of this paper. These ideas are generated by
the numerical experience of the author in working with the convexification
for a CIP with single measurement data \cite{KNT} as well as for an
ill-posed Cauchy problem for a quasilinear parabolic PDE \cite{KB}.

First, even though the above theory is valid only for sufficiently large
values of $\lambda ,$ in fact, $\lambda \in \left[ 1,3\right] $ worked well
in \cite{KB,KNT}. Another observation is that it is more effective to work
with such functions $\xi \left( x\right) $, which are simple and change
rather slowly. However, the function $\xi \left( x\right) $ in (\ref{6.4})
changes rapidly due to the presence of the parameter $\nu >1.$ On the other
hand, it was heuristically established in \cite{Liem,Liem2} that the
stability of the numerical solution of an analog of CISP1 (section 5.2) can
be improved if the Dirichlet and Neumann data at the backscattering side $%
\Gamma _{b}$ of $\Omega $ are complemented on the rest of $\partial \Omega $
by the Dirichlet data generated by the solution of the problem (\ref{6.5}), (%
\ref{6.6}) for the case $c\left( x\right) \equiv 1.$ At the same time, it
was also observed in \cite{Liem,Liem2} that this complement does influences
the accuracy of the solution insignificantly.

One can prove that the CWF in the latter case can be chosen as

$\varphi _{\lambda }^{\left( 1\right) }\left( x\right) =\exp \left[ \lambda
\left( x_{3}-B-b_{1}\right) ^{2}\right] ,$ where $b_{1}>0$ is any number.
One can simplify even this choice via choosing another CWF as $\varphi
_{\lambda }^{\left( 2\right) }\left( x\right) =\exp \left( -\lambda
x_{3}\right) :$ this CWF works for the 1-D operator $d^{2}/dx_{3}^{2},$ see
lemma 6.1 in \cite{KNT}. Then, however, one needs to assume that all
derivatives with respect to $x_{1}$ and $x_{2}$ are written in finite
differences, unlike derivatives with respect to $x_{3}.$ In this case, the
parameter $\lambda $ would depend on the grid step size. One can proceed
similarly in the case of CISP2.

Assume now that the restricted DN data for the elliptic case are given on
the sphere $S=\left\{ \left\vert x\right\vert =1\right\} $ and that $\Omega
=\left\{ \kappa <\left\vert x\right\vert <1\right\} ,$ where $\kappa
=const.\in \left( 0,1\right) .$\ In addition, assume that $A_{0}=\Delta .$
Writing the operator $\Delta $ in spherical coordinates as $\Delta
_{r,\varphi ,\theta },$ where $r\in \left( \kappa ,1\right) ,\varphi \in
\left( 0,2\pi \right) ,\theta \in \left( 0,\pi \right) ,$ one can prove an
analog of the Carleman estimate (\ref{2.90})-(\ref{2.92}) for $\left( \sqrt{r%
}\Delta _{r,\varphi ,\theta }u-2u_{r}/\sqrt{r}\right) ^{2}e^{2\lambda r}.$
However, when integrating the analog of (\ref{2.90}) over $\Omega ,$ one
should replace the conventional $r^{2}\sin \theta drd\varphi d\theta $ with $%
\sin \theta drd\varphi d\theta .$ As to the term $\left( -2u_{r}/\sqrt{r}%
\right) ,$ recall that Carleman estimates are independent on terms with
derivatives, whose order is less than the order of derivatives in the
principal part of a corresponding PDE operator, see page 39 in \cite{KT}.
Hence, the function $\varphi _{\lambda }^{\left( 3\right) }\left( r\right)
=e^{2\lambda r}$ might be an appropriate choice of the CWF in this case. The
2-D case is similar. It is worthy to repeat now that if the data are given
rather far from the domain $\Omega ,$ then the data propagation procedure is
recommended, see \cite{Liem2} for a detailed description.

Considerations about the CWF, which are similar with the ones above in this
section, can be also brought in for the case of CIPs with restricted DN data
for parabolic PDEs (section 6). Here is an example: Let $n=3$ and let the
domain $\Omega _{1}$ be the same as the domain $\Omega $ in (\ref{6.60}).
Also, assume that the restricted DN data (\ref{7.7}) are given at $\left(
\partial \Omega _{1}\cap \left\{ x_{3}=0\right\} \right) \times \left( \zeta
,T\right) =\Psi \times \left( \zeta ,T\right) .$ In addition, assume that we
have Dirichlet data at $\left( \partial \Omega _{1}\diagdown \Psi \right)
\times \left( \zeta ,T\right) =\Phi $ and that in (\ref{7.2}) the operator $%
A_{0}=\partial _{t}-\Delta $. In this case the CWF of section 6 can be
simplified as $\varphi _{\lambda }\left( x,t\right) =\exp \left[ \lambda
\left( \left( x_{3}-B-b_{1}\right) ^{2}-\left( t-T/2\right) ^{2}\right) %
\right] ,$ where $b_{1}>0$ is any number.

In some applications the point source might run along a circle surrounding
either the entire domain $\Omega $ in the 2-D case or a 2-D cross-section of 
$\Omega $ in the 3-D case. The above process can be modified as follows
then: Let $\left\{ \left( r,\varphi ,\theta \right) :r=1,\varphi \in \left(
0,2\pi \right) ,\theta =\pi /2\right\} $ be that circle. Choose a small
number $\varepsilon \in \left( 0,\pi \right) .$ Next, in the process of
sections 2-4, replace $x_{0}\in \left[ 0,1\right] $ in (\ref{2.14}) with $%
\varphi _{0}\in \left[ \varepsilon ,2\pi -\varepsilon \right] $ and modify
the orthonormal basis of section 2.3 accordingly.

\begin{center}
\textbf{Acknowledgments}
\end{center}

This work was supported by US Army Research Laboratory and US Army Research
Office grant W911NF-15-1-0233 and by the Office of Naval Research grant
N00014-15-1-2330.

\end{document}